\documentclass[final,a4paper,reqno]{jmsj}
\usepackage[pdflatex,dvipsnames]{xcolor}
\usepackage{mathrsfs,amsmath,amssymb,amsthm}
\usepackage[dvipdfmx]{graphicx}
\usepackage{ulem}
\usepackage[]{showlabels}

\newcommand{\R}{\mathbb R}

\theoremstyle{thmstyleone}%
\newtheorem{theorem}{Theorem}
\theoremstyle{thmstyletwo}%
\newtheorem{remark}{Remark}%

\theoremstyle{thmstylethree}%
\newtheorem{definition}{Definition}%

\newtheorem{lemma}[theorem]{Lemma}

\begin{document}

\title[Hadamard variation of eigenvalues]
{Hadamard variation of eigenvalues with
 respect to general domain perturbations}

\author{Takashi \textsc{Suzuki}}% Author Name (\sc should NOT be used here)
\address{Center of Mathematical Modeling and Data Science, \\
Osaka University, \\
Toyonaka, 560-8531, Japan}
\email{suzuki@sigmath.es.osaka-u.ac.jp}

\author{Takuya \textsc{Tsuchiya}}% Author Name (\sc should NOT be used here)
%\address{Graduate School of Science and Engineering, \\
%Ehime University, \\
%Matsuyama 790-8577, Japan}
\address{Center of Mathematical Modeling and Data Science, \\
Osaka University, \\
Toyonaka, 560-8531, Japan}
\email{tsuchiya.takuya.plateau@kyudai.jp}

% Subject code(s)

%\recdate{20XX}{11}{99}% Date of reception
%\revdate{20XX}{1}{99}%  Date of revision

\begin{abstract}
 We study Hadamard variation of eigenvalues of Laplacian with 
respect to general domain perturbations. We show their existence up to the 
second order rigorously and characterize the derivatives, using
associated eigenvalue problems in  finite dimensional spaces. Then smooth rearrangement of multiple eigenvalues is explicitly given. This
result follows from an abstract theory, applicable to general
perturbations of symmetric bilinear forms.
\end{abstract}

\keywords{eigenvalue problem, perturbation theory of linear operators, domain deformation, Hadamard's
variational formulae, Garabedian-Schiffer's formula\\
\hspace{4mm}
This work was promoted in RIMS program for joint
research during 2019-2021. The  authors thank Professors Hideyuki Azegami
and Erika Ushikoshi for many detailed discussions at these occasions. The first
author was supported by JSPS Grant-in-Aid for Scientific Research
19H01799. The second author was supported by JSPS Grant-in-Aid for
Scientific Research 21K03372}

\maketitle

\section{Introduction}\label{sec1}  

Our purpose is to study Hadamard variation of eigenvalues of the Laplace operator with mixed boundary conditions. We characterize the first and the second derivatives, using associated finite dimensional eigenvalue problems, particularly, for multiple eigenvalues. 

Hadamard's variational formulae are used to provide effective numerical schemes for shape optimization and free boundary problems. In \cite{st05, st11} we have introduced a new variational formulation to filtration problem and applied the Hadamard variation. We have also derived Hadamard's variational formulae for Green's functions of the Poisson equation, extending Liouville's volume and area formulae up to the second order \cite{SuzTsu16, st22}.  This paper is devoted to the eigenvalue problem concerning general perturbation of Lipschitz domains. Meanwhile we develop abstract theory applicable to other problems. 

So far, $C^1$ and analytic rearrangements of multiple eigenvalues of
self-adjoint operator have been discussed, for example, pp. 44--52 of
\cite{rellich53} and p.487 of \cite{ch82}. There are, however, several
comments on the complexity of the proof  on $C^1$
category, as in pp. 122--123 of \cite{kato76} and p. 490 of \cite{ch82}.  In this paper we  examine  this process of rearrangement in details, using the above described characterization of the derivatives, to ensure their efficiencies in $C^2$ categories.

Let $\Omega$ be a bounded Lipschitz domain in $n$-dimensional Euclidean space $\R^n$ for $n \ge 2$. Suppose that its boundary $\partial\Omega$ is divided into two relatively open  
sets $\gamma_0$ and $\gamma_1$, satisfying 
\begin{equation}
   \overline{\gamma_0} \cup \overline{\gamma_1}= \partial\Omega,
    \quad
   \overline{\gamma_0} \cap \overline{\gamma_1} = \emptyset.
  \label{assum-gamma}
\end{equation}
As in \cite{st22}, we thus avoid the technical difficulty described in p. 272 of \cite{fss01}, when $\gamma_i$, $i=0,1$ have their boundaries on $\partial\Omega$.  

We study the eigenvalue problem of the Poisson equation with mixed boundary condition, 
\begin{equation}
  - \Delta u = \lambda u \ \mbox{in $\Omega$}, \quad
      u = 0 \ \mbox{on $\gamma_0$}, \quad
      \frac{\partial u}{\partial \nu} = 0 \ \mbox{on $\gamma_1$},
   \label{2}
\end{equation}
where
\[ \Delta=\sum_{i=1}^n \frac{\partial^2}{\partial x_i^2} \] 
is the Laplacian and $\nu$ denotes the outer unit normal vector on $\partial\Omega$. This problem takes the weak form, finding $u$ satisfying  
\begin{equation} 
u\in V, \ B(u,u)=1, \quad A(u,v)=\lambda B(u,v), \ \forall v\in V 
 \label{3}
\end{equation} 
defined for  
\[ A(u,v)=\int_\Omega \nabla u\cdot \nabla v \ dx, \quad B(u,v)=\int_\Omega uv \ dx \] 
and 
\begin{equation} 
V=\{ v\in H^1(\Omega) \mid \left. v\right\vert_{\gamma_0}=0\}.   
 \label{v}
\end{equation} 
This $V$ is a closed subspace of $H^1(\Omega)$ under the norm 
\[ \Vert v\Vert_V=\sqrt{\Vert \nabla v\Vert_2^2+\Vert v\Vert_2^2}. \]  

To justify the above weak formulation, we confirm several fundamental facts on the Lipschitz domain \cite{kjf77}. First, the trace operator 
\[ v\in C^\infty(\overline{\Omega})\mapsto \left. v\right\vert_{\partial\Omega} \in C(\partial\Omega) \] 
defined for 
\[ C^\infty(\overline{\Omega})=\{ v\in \overline{\Omega} \rightarrow \R \mid \exists \tilde v\in C_0^\infty(\R^n), \ \left. \tilde v\right\vert_{\overline{\Omega}}=v \} \] 
is so extended as a bounded linear operator    
\[ v\in H^1(\Omega) \ \mapsto \ \left. v\right\vert_{\partial\Omega}\in H^{1/2}(\partial\Omega), \]  
and there arises the isomorphism 
\[ v\in V/H^1_0(\Omega) \ \mapsto \left. v\right\vert_{\gamma_1}\in H^{1/2}(\gamma_1). \] 
Second, if $\Delta v\in V'$ is satisfied in the sense of distributions in $\Omega$, the normal derivative of $v\in V$ on $\gamma_1$ is defined as in   
\[ \frac{\partial v}{\partial \nu}\in H^{-1/2}(\gamma_1)=(H^{1/2}(\gamma_1))', \] 
and it holds that   
\[ \left\langle \varphi, \frac{\partial v}{\partial \nu}\right\rangle_{H^{1/2}(\gamma_1), H^{-1/2}(\gamma_1)}=(\nabla v, \nabla\varphi)_{L^2(\Omega)}+\langle \varphi, \Delta v\rangle_{V, V'}, \quad \forall \varphi\in V,  \] 
where $\langle \cdot, \cdot \rangle_{Y,Y'}$ denotes the paring between the Banach space $Y$ and its dual space $Y'$. See Theorem 2 of \cite{st22}. 

To confirm the well-posedness of (\ref{3}), we note, first, that 
if $\gamma_0\neq\emptyset$, there is coercivity of $A:V\times V\rightarrow \R$, which means the existence of $\delta>0$ such that  
\begin{equation} 
A(v,v)\geq \delta \Vert v\Vert_V^2, \quad \forall v\in V. 
 \label{4}
\end{equation} 
If $\gamma_0=\emptyset$ we replace $A$ by $A+B$, denoted by $\tilde A$. Then this $\tilde A:V\times V\rightarrow \R$ is coercive, and the eigenvalue problem 
\[ u\in V, \ B(u,u)=1, \quad \tilde A(u,v)=\tilde \lambda B(u,v), \ \forall v\in V, \] 
is equivalent to (\ref{3}) by $\tilde \lambda=\lambda+1$.  Henceforth, we assume (\ref{4}), using this  reduction if it is necessary. 

Second, we note that $A:V\times V\rightarrow {\bf R}$ and $B:X\times X\rightarrow {\bf R}$ are bounded, coercive, and symmetric bilinear forms for $X=L^2(\Omega)$. Since $V\hookrightarrow X$ is compact, there is a sequence of eigenvalues to (\ref{3}), denoted by 
\[ 0<\lambda_1\leq \lambda_2\leq \cdots \rightarrow +\infty. \] 
The associated eigenfunctions, $u_1, u_2, \cdots$, furthermore, form a complete ortho-normal system in $X$, provided with the inner product induced by $B=B(\cdot, \cdot)$: 
\[ B(u_i, u_j)=\delta_{ij}, \quad A(u_j, v)=\lambda_jB(u_j, v), \ \forall v \in V, \quad i,j=1,2,\cdots. \]  
The $j$-th eigenvalue of (\ref{3}) is given by the mini-max principle 
\begin{equation} 
\lambda_j=\min_{L_j}\max_{v\in L_j\setminus \{0\}}R[v]=\max_{W_j}\min_{v\in W_j\setminus \{0\}}R[v], 
 \label{minimax}
\end{equation} 
where 
\[ R[v]=\frac{A(v,v)}{B(v,v)} \] 
is the Rayleigh quotient, and $\{ L_j\}$ and $\{ W_j\}$ denote the families of all subspaces of $V$ with dimension and codimension $j$ and $j-1$, respectively.  See  \cite{suzuki22}, for example, for these fundamental facts.    

Let  
\begin{equation} 
T_t : \Omega\rightarrow \Omega_t= T_t(\Omega), \quad \vert t\vert < \varepsilon_0   
 \label{13}
\end{equation} 
be a family of bi-Lipschitz homeomorphisms for $\varepsilon_0>0$, satisfying $T_0=I$, the identity mapping. We assume that $T_tx$ is continuous in $t$ uniformly in $x\in \Omega$, and recall the following definition used in \cite{st22}.  

\begin{definition}\label{def1} 
The family  $\{ T_t\}$ of bi-Lipschitz homeomorphisms is said to be $p$-differentiable in $t$ for $p\geq 1$, if $T_tx$ is  $p$-times differentiable in $t$ for any $x\in\Omega$ and the mappings 
\[ 
\frac{\partial^\ell}{\partial t^\ell}DT_t, \ \frac{\partial^\ell}{\partial t^\ell}(DT_t)^{-1}:\Omega \rightarrow M_n(\R), \quad 0\leq \ell \leq p  \nonumber\\  
\] 
are uniformly bounded in $(x,t)\in \Omega\times (-\varepsilon_0, \varepsilon_0)$, where $DT_t$ denotes the Jacobi matrix of $T_t:\Omega \rightarrow \Omega_t$ and $M_n(\R)$ stands for the set of real $n\times n$ matrices. This $\{ T_t\}$ is said to be continuously $p$-differentiable in $t$ if it is $p$-differentiable and the mappings 
\[ 
t\in (-\varepsilon_0, \varepsilon_0)\mapsto \frac{\partial^\ell}{\partial t^\ell}DT_t, \ \frac{\partial^\ell}{\partial t^\ell}(DT_t)^{-1} \in L^\infty(\Omega \rightarrow M_n(\R)), \quad 0\leq \ell \leq p  \] 
are continuous.  
\end{definition}

Putting  
\begin{equation} 
T_t(\gamma_i)=\gamma_{it}, \quad i=0,1, 
 \label{6}
\end{equation} 
we introduce the other eigenvalue problem 
\begin{equation}
  - \Delta u= \lambda u \ \mbox{in $\Omega_t$}, \quad
      u = 0 \ \mbox{on $\gamma_{0t}$}, \quad
      \frac{\partial u}{\partial \nu} = 0 \ \mbox{on $\gamma_{1t}$},
   \label{7}
\end{equation}
which is reduced to finding 
\begin{equation} 
u\in V_t, \ \int_{\Omega_t}u^2 \ dx=1, \quad \int_{\Omega_t}\nabla u\cdot \nabla v \ dx=\lambda\int_{\Omega_t}uv \ dx, \ \forall v\in V_t 
 \label{88}
\end{equation} 
for 
\begin{equation} 
V_t=\{ v\in H^1(\Omega_t) \mid \left. v\right\vert_{\gamma_{0t}}=0\}.  
 \label{9}
\end{equation} 
Let $\lambda_j(t)$ be the $j$-th eigenvalue of the eigenvalue problem (\ref{7}). In Section \ref{sec2} we confirm that the eigenvalue problem (\ref{88})-(\ref{9}) is reduced to 
\begin{equation} 
u\in V, \ B_t(u,u)=1, \quad A_t(u,v)=\lambda B_t(u,v), \ \forall v\in V 
 \label{wf}
\end{equation} 
by the transformation of variables $y=T_tx$ for $V\subset H^1(\Omega)$ defined by (\ref{v}), where 
\begin{equation} 
B_t(u,v)=\int_\Omega uv a_t \ dx, \quad A_t(u,v)=\int_\Omega Q_t[\nabla u, \nabla v] a_t \ dx,  
 \label{ab}
\end{equation} 
and 
\begin{equation} 
a_t=\det DT_t, \quad Q_t=(DT_t)^{-1}(DT_t)^{-1T}. 
 \label{cd}
\end{equation}

Several representation formulae for 
\[ \lambda'_j(t)=\lim_{h\rightarrow 0}\frac{1}{h}(\lambda_j(t+h)-\lambda_j(t)) \] 
and 
\[ \lambda_j''(t)=\lim_{h\rightarrow 0}\frac{1}{h}(\lambda'_j(t+h)-\lambda'_j(t)) \] 
have been derived in \cite{gs52}. Here we characterize these derivatives, using associated finite dimensional eigenvalue problems (Theorems \ref{thm2.1} and \ref{thmsecond}). We prove also the existence of derivatives, particularly, the second derivatives of multiple eigenvalues. These results follow from the differentiability of $A_t$ and $B_t$ in $t$, defined by 
\begin{equation} 
\dot A_t(u,v)=\frac{d}{dt}A_t(u,v),  \ \dot B_t(u,v)=\frac{d}{dt}B_t(u,v), \quad u, v\in V 
 \label{dot}
\end{equation} 
and 
\begin{equation} 
\ddot A_t(u,v)=\frac{d^2}{dt^2}A_t(u,v),  \ \ddot B_t(u,v)=\frac{d^2}{dt^2}B_t(u,v), \quad u, v\in V.    
 \label{ddot+}
\end{equation}  
Then we show $C^1$ and $C^2$ rearrangements of eigenvalues, if  these bilinear forms are continuous in $t$ (Theorems \ref{rellich+} and \ref{c2theorem}). We give the algorithm explicitly, as the {\it transversal rearrangement} in Definition \ref{transversal}. Consequently, no rearrangment is necessary to confirm $C^1$ or $C^2$ smoothness of eigenvalues, if their multiplicities are constant. Also elementary symmetric functions made by possible multiple eigenvalues are $C^1$ or $C^2$. These properties are noticed by \cite{pdm1} and \cite{pdm2} in the real analytic category.

In \cite{ju15}, unilateral derivatives of the first order of the eigenvalue of the Stokes operator are calculated. It  characterizes the derivatives using boundary integrals when the deformation of the domain is of normal direction. Our abstract theory reproduces and extends this result by the use of Piola transformation described in \cite{mh94}. 

Our argument is executed in the $H^1$ category without requiring any further elliptic regularities. Hence the Lipschitz continuity of $\partial\Omega$ is sufficient to ensure all the results on (\ref{7}) under the general perturbation of domains. We recall, in this context, that numerical computations on partial differential equations  are mostly executed on Lipschitz domains. There, perturbation of the domain using Lipschitz continuous vector fields is often applied. This is the method of trial domains, efficient even to the case that {\it normal perturbation} of domains, called in \cite{st05}, does not work because of the presence of corners on the boundary \cite{st05, st11}. 

\section{Summary}\label{summary0}

As is noted in the previous section, $C^1$ and analytic categories for the smoothness of $\lambda_j(t)$ in $t$ have been discussed. Since each  eigenvalue $\lambda$ of (\ref{wf}) is isolated, we can reduce this problem to a finite dimensional eigenvalue problem by the Lyapunov-Schmidt reduction as in p.486 of \cite{ch82}.  Rellich \cite{rellich53} in p.45 showed in this case that the continuous differentiablity in $t$ of $A_t$ and $B_t$ implies that of $\lambda_j(t)$ under a suitable change of its order. Kato \cite{kato76} in p.123 then provided an alternative proof of this $C^1$ rearrangement. 

The other category of analyticity is studied in Chapter II of \cite{rellich53} and p.370 of \cite{kato76}.  If $A_t$ and $B_t$ are analytic in $t$, the eigenvlue problem is reduced to an algebraic equation with analytic coefficients in $t$. The Puiseux expansion of $\lambda_j(t)$ at $t=0$ is described in p.31 of \cite{rellich53} and Chapter 2, Section 1 of \cite{kato76}. Hence this $\lambda_j(t)$ is realized as an analytic function defined on a Riemann surface. 

We study $C^2$ category. To begin with, we confirm Rellich's theorem on $C^1$ category, the continuous differentiablity of rearranged eigenvalues. Here we show this rearrangement explicitly, to reach existence, characterization, and continuity of the second derivatives (Definition \ref{transversal}). In more details, first, the existence of $\dot A_t$ and $\dot B_t$ in (\ref{dot}) implies that of the first unilateral derivatives  
\begin{equation} 
\dot \lambda^\pm_j(t)=\lim_{h\rightarrow \pm 0}\frac{1}{h}(\lambda_j(t+h)-\lambda_j(t)) 
 \label{unilateral}
\end{equation} 
for each $j$ and $t$. These derivatives, furthermore, are characterized by the other finite dimensional eigenvalue problems in accordance with the multiplicity of $\lambda_j(t)$ (Theorem \ref{thm4}).
Second, if the above $\dot A_t$ and $\dot B_t$ are continuous in t, and if 
\begin{equation} 
\lambda_{k-1}(t)<\lambda_k(t)\leq \cdots \leq \lambda_{k+m-1}(t)<\lambda_{k+m}(t) 
 \label{x19}
\end{equation} 
holds for $t\in I= (-\varepsilon_0, \varepsilon_0)$, there are $C^1$ curves  
\[ \tilde C_j, \ k\leq j\leq k+m-1, \] 
made by at most countably many rearrangements of the $C^0$ curves  
\[ C_j=\{ \lambda_j(t) \mid t\in I\}, \quad k\leq j\leq k+m-1. \] 
 (Theorem \ref{c1theorem}).  

In this paper, we notice two properties of the unilateral derivatives $\dot \lambda_j^\pm$, for the proof of this Rellich's theorem. First, there arises that 
\[ \dot \lambda_j^+(t)=\dot \lambda_{2\ell+n-1-j}^-(t), \quad \ell\leq j\leq \ell+n-1 \] 
if 
\begin{equation} 
\lambda_{\ell-1}(t)<\lambda_\ell(t)=\cdots = \lambda_{\ell+n-1}(t)<\lambda_{\ell+n}(t)   
 \label{18}
\end{equation} 
holds for $\ell \geq k$, $\ell+n\leq m$ (Theorem \ref{thm4}). Second, the unilateral derivative $\dot \lambda_j^\pm$ are provided with the unilateral continuity as in  
\[ \lim_{h\rightarrow \pm 0}\dot \lambda_j^\pm(t+h)=\dot \lambda_j^\pm(t) \] 
if both $\dot A_t$ and $\dot B_t$ are continuous in $t\in I$ (Theorem \ref{abstconti}). 

As for the second derivatives of eigenvalues, we assume the existence of $\ddot A_t$ and $\ddot B_t$ in (\ref{ddot+}) besides $\dot A_t$ and $\dot B_t$. Then each $j=1,2,\cdots$ admits the existence of
\begin{equation}    
\ddot \lambda_j^\pm(t)=\lim_{h\rightarrow \pm 0}\frac{2}{h^2}(\lambda_j(t+h)-\lambda_j(t)-h\dot \lambda_j^\pm(t)). 
 \label{18+}
\end{equation} 
These limits are again unilateral and characterized by the other eigenvalue problem in a finite dimensional space (Theorem \ref{secondabst}). The unilateral continuity of these $\ddot \lambda_j^\pm(t)$ is then assured under the continuity of $\ddot A_t$ and $\ddot B_t$ in $t$, similarly. These properties induce $C^2$ smoothness of $\tilde C_j$, $k\leq j\leq k+m-1$, once their $C^1$ smoothness is achieved (Theorem \ref{finaltheorem}).

This paper is composed of eight sections. The results on the Hadamard variation to (\ref{7}) are described in the next section. In Section \ref{sec2} we use the transformation of variables to introduce an abstract setting of the problem. Section \ref{sec3} is concerned on the continuity of eigenvalues and eigenspaces. We study the first derivative of eigenvalues in Section \ref{sec5}, and then its $C^1$ rearrangement in Section \ref{sec:rearrange}.  Finally, Section \ref{sec6} is devoted to the second derivatives.

\section{Hadamard variation}\label{summary} 

Here we state the results on the first and the second derivatives of $\lambda_j=\lambda_j(t)$ in (\ref{wf}).  Fix $t\in I$, and assume (\ref{18}) for $\ell=k$ and $n=m$ for $k,m=1,2,\cdots$ with  the convention $\lambda_{0}(t)=-\infty$. Put   
\begin{equation}
\lambda\equiv \lambda_k(t)=\cdots=\lambda_{k+m-1}(t),  
 \label{put}
\end{equation} 
and let $Y^\lambda_t$, $\dim Y^\lambda_t=m$, be the eigenspace corresponding to this eigenvalue $\lambda$. 

\begin{theorem}\label{thm2.1} 
Assume the above situation, and let $\{T_{t'}\}$ be $1$-differentiable at $t'=t$. Then, there exist the unilateral limits  
\[ \dot\lambda_j^\pm=\lim_{h\rightarrow \pm 0}\frac{1}{h}(\lambda_j(t+h)-\lambda), \quad k\leq j\leq k+m-1, \] 
which satisfies  
\begin{equation} 
\nu_j\equiv \dot \lambda_j^+=\dot \lambda_{2k+m-1-j}^-, \quad k\leq j\leq k+m-1. 
 \label{derivative}
\end{equation} 
This $\nu_j$ is the $q$-th eigenvalue of the matrix 
\begin{equation} 
G^\lambda_t=\left( E^\lambda_t(\tilde \phi_i, \tilde \phi_j)\right)_{1\leq i, j\leq m} 
 \label{glambda}
\end{equation} 
for $q=j-k+1$, where $\{ \tilde \phi_j \mid 1\leq j\leq m\}$ is a basis of $Y^\lambda_t$, satisfying 
\begin{equation} 
B_t(\tilde \phi_i, \tilde \phi_j) =\delta_{ij}, \quad 1\leq i, j\leq m 
 \label{deltil}
\end{equation} 
and 
\[ E^\lambda_t=\dot A_t-\lambda \dot B_t \]
for $\dot A_t$ and $\dot B_t$ defined by (\ref{ab}), (\ref{cd}), and (\ref{dot}).  
\end{theorem} 

\begin{remark}\label{orthgonalt}
If $\langle \phi_j \mid 1\leq j\leq m\rangle$ is the other basis of $Y^\lambda_t$ satisfying 
\begin{equation}  
B_t(\phi_i, \phi_j)=\delta_{ij}, \quad 1\leq i, j\leq m, 
 \label{xx1}
\end{equation} 
it holds that 
\[ \phi_j=\sum_{i=1}^mq^i_j\tilde \phi_i, \ 1\leq j\leq m \] 
with the orthogonal $m\times m$ matrix $Q=(q^i_j)$. Hence $\nu_j$, $k\leq j\leq k+m-1$, in (\ref{derivative}) is determined, indpendent of the choice of $\langle \tilde \phi_j \mid 1\leq j\leq m\rangle$.  
\end{remark}

\begin{remark} 
Under (\ref{put}), it holds that 
\[ \dot \lambda_k^+(t)\leq \cdots \leq \dot \lambda_{k+m-1}^+(t), \quad \dot \lambda_k^-(t)\geq \cdots \geq \dot \lambda_{k+m-1}^-(t). \] 
\end{remark}

A direct consequnece of this theorem is the existence of the unilateral derivatives 
$\dot\lambda_j^\pm(t)$ in (\ref{unilateral}) for any $t\in I$ and $j=1,2,\cdots$, if $\{ T_t\}$ is $1$-differentiable in $I$. Then we obtain the following theorem. 

\begin{theorem}\label{thm2+}
If $\{ T_t\}$ is continuously $1$-differentiable in $t\in I$, the functions $\dot\lambda^\pm_j=\dot\lambda^\pm_j(t)$ are unilaterally continuous, so that it holds that 
\[ \lim_{h\rightarrow \pm 0}\dot \lambda_j^\pm(t+h)=\dot \lambda_j^\pm(t) \] 
for any $t$ and $j$.  
\end{theorem} 

This fact ensures the following theorem of Rellich. 

\begin{theorem}\label{rellich+} 
Let $\{ T_t\}$ be continuously $1$-differentiable in $t$, and assume 
\[ \lambda_{k-1}(t)<\lambda_k(t)\leq \cdots \leq \lambda_{k+m-1}(t)<\lambda_{k+m}(t), \quad \forall t\in I  \] 
for some $k,m=1,2,\cdots$. Let 
\[ C_j=\{ \lambda_j(t) \mid t\in I\}, \quad k\leq j\leq k+m-1 \] 
be $C^0$ curves. Then, there exist $C^1$ curves denoted by $\tilde C_i$, $k\leq i\leq k+m-1$, made by a rearrangement of 
\[ \{ C_j\mid k\leq j\leq k+m-1\} \] 
at most countably many times.  
\end{theorem} 

Turning to the second derivatives, we fix $t\in I$ again, and assume that $\{ T_{t'}\}$ is twice differentiable at $t'=t$. Suppose (\ref{18}) for $\ell=k$ and $n=m$, put $\lambda$ as in (\ref{put}), and let $k\leq \ell <r\leq k+m$ be such that 
\begin{equation} 
\dot\lambda^+_{\ell-1}<\lambda'\equiv \dot \lambda^+_\ell=\cdots =\dot \lambda^+_{r-1}<\dot \lambda^+_r   
 \label{eq15+}
\end{equation} 
in Theorem \ref{thm2.1}. To state the finite dimensonal eigenvalue problem characterizing the second derivatives of $\lambda_j(t')$ at $t'=t$ for $\ell \leq j\leq r-1$, we introduce the following definition. 

\begin{definition}\label{def2} 
Let $R:X=L^2(\Omega)\rightarrow Y^\lambda_t$ be the orthogonal projection with respect to $B_t(\cdot, \cdot)$ and $P=I-R$, where $I:X\rightarrow X$ is the identity operator. Let, furthermore, $A_t$, $B_t$, $\dot A_t$, $\dot B_t$, $\ddot A_t$, and $\ddot B_t$ be as in (\ref{ab}), (\ref{cd}), (\ref{dot}), and (\ref{ddot+}).  
Then we define $w=\gamma(u)\in PV$ for $u\in V$ by 
\begin{equation} 
C_t(w,v)=-\dot C_t^{\lambda, \lambda'}(u, v), \quad \forall v\in PV,  
 \label{wc0}
\end{equation} 
where 
\[ C_t=A_t-\lambda B_t, \quad \dot C_t^{\lambda, \lambda'}=\dot A_t-\lambda \dot B_t-\lambda'B_t.  \] 
We put also   
\[ F_t^{\lambda, \lambda'}(u,v)=(\ddot A_t-\lambda \ddot B_t-2\lambda'\dot B_t)(u,v)-2C_t(\gamma(u), \gamma(v)), \quad u, v\in V. \]  
\end{definition} 

\begin{remark}\label{remgamma} 
To confirm the unique solvability of $w=\gamma(u)$, let $Q:L^2(\Omega)\rightarrow Z^k_t$ be the orthogonal projection with respect to $B_t(\cdot,\cdot)$, where $Z_t^k$ denotes the finite dimensional space of $L^2(\Omega)$ generated by the eigenfunctions corresponding to the eigenvalues $\lambda_1(t), \cdots, \lambda_{k-1}(t)$. Let, furthermore, $V_0=QV$ and $V_1=(I-Q)RV$. First, there is a unique $w_0\in V_0$ satisfying 
\begin{equation} 
C_t(w_0,v)=-\dot C_t^{\lambda, \lambda'}(u,v), \quad \forall v\in V_0 
 \label{wc1}
\end{equation} 
because $C_t$ is negative definite on $V_0\times V_0$. Second, there is also a unique $w_1\in V_1$ satisfying 
\begin{equation}  
C_t(w_1, v)=-\dot C_t^{\lambda, \lambda'}(u,v), \quad \forall v\in V_1  
 \label{wc2}
\end{equation} 
because $C_t$ is positive definite on $V_1\times V_1$. Then we obtain (\ref{wc0}) for $w=w_0+w_1$, because  
\[ C_t(w_0, v)=C_t(w_0, Qv), \ C_t(w_1,v)=C_t(w_1, (I-Q)v), \quad \forall v\in RV, \] 
and hence 
\begin{eqnarray*} 
C_t(w, v) & = & C_t(w_0, v)+C_t(w_1, v)= C_t(w_0, Qv)+C_t(w_1, (I-Q)v) \\ 
& = & -\dot C_t^{\lambda, \lambda'}(u, Qv)-\dot C_t^{\lambda, \lambda'}(u,(I-Q)v)=-\dot C_t(u,v), \quad \forall v\in RV. 
\end{eqnarray*}  
We thus obtain $w\in RV$ satisfying (\ref{wc0}). If $w\in RV$ is a solution to (\ref{wc0}), conversely, then $w_0=Qw\in V_0$ and $w_1=(I-Q)w\in V_1$ solve (\ref{wc1}) and (\ref{wc2}), respectively. Then there arises the uniqueness of such $w=w_0+w_1$ because these $w_0\in V_0$ and $w_1\in V_1$ are unique. 
\end{remark}  

Recall 
\[ Y^\lambda_t=\langle \tilde\phi_j \mid k\leq j\leq k+m-1 \rangle \] 
with (\ref{deltil}).  

\begin{theorem}\label{thmsecond}
Under the above situation of (\ref{18}) and (\ref{eq15+}), there exist  
\begin{equation} 
\lambda''_j=\lim_{h\rightarrow +0}\frac{2}{h^2}(\lambda_j(t+h)-\lambda-h\lambda'), \quad \ell \leq j\leq r-1.  
 \label{2deri}
\end{equation}
This $\lambda''_j$ is the $q$-th eigenvalue of the matrix 
\begin{equation} 
H^{\lambda, \lambda'}_t=\left( F^{\lambda, \lambda'}_t(\tilde \phi_i, \tilde \phi_j)\right)_{\ell \leq i, j\leq r-1} 
 \label{hlambda}
\end{equation} 
for $q=j-\ell+1$. 
If (\ref{eq15+}) is replaced by 
\[ \dot\lambda^-_{\ell-1}(t)>\lambda'\equiv \dot \lambda^-_{\ell}(t)=\cdots=\dot \lambda_{r-1}^-(t)>\dot \lambda_{r}^-(t),  \] 
there arises that   
\[ \lambda_j''(t)=\lim_{h\rightarrow -0}\frac{2}{h^2}(\lambda_j(t+h)-\lambda-h\lambda'), \quad \ell \leq j\leq r-1.  \] 
\end{theorem} 

As in Remark \ref{orthgonalt} on the first derivative, the above $\lambda''_j$, $\ell\leq j\leq r-1$, are determined  independent of the choice of 
\[ \langle \tilde \phi_j \mid \ell \leq j\leq r-1\rangle. \] 

Theorems \ref{thm2.1}, \ref{thmsecond} imply also the existence of the unilateral limits $\ddot \lambda_j^\pm(t)$ in (\ref{18+}) for any $t$ and $j$ if $\{ T_t\}$ is $2$-differentiable in $t\in I$, that is, 
\[ \ddot \lambda_j^\pm(t)=\lim_{h\rightarrow \pm 0}\frac{2}{h^2}(\lambda_j(t+h)-\lambda_j(t)-h\dot \lambda_j^\pm(t)). \] 

\begin{remark} 
By Liouville's theorem on general perturbation of domains studied in
 \cite{st22}, the matrix  $G_t^\lambda$ in (\ref{glambda}) is
 represented by the surface integrals of $\tilde \phi_j$,  $1\leq j\leq
 m$. This property is confirmed by \cite{ju15} for the Stokes operator
 with Dirichlet condition under a special perturbation of domains,
 called the {\it normal perturbation} in \cite{SuzTsu16}. Similarly, the
 matrix $H_t^{\lambda, \lambda'}$ in (\ref{hlambda}) is represented by
 the surface integrals of $\tilde \phi_j$ and $\gamma(\tilde \phi_j)$,
 $1\leq j\leq m$. %In particular, the first and the second derivatives of
 %the first eigenvalue, $\lambda_1$, are represented by the surface
 %integrals of $\tilde \phi_1$, because in this case there arises that $Z_t^k=\{0\}$ in Remark \ref{remgamma}, besides the simplicity of $\lambda_1$. This fact on the first eigenvalue is noticed by  \cite{gs52} for (\ref{2}) with $\gamma_1=\emptyset$ in two space dimension under the first order perturbation of the domain in $t$. 
\end{remark}

Then we obtain the following theorems. 

\begin{theorem}\label{secondconti}
If $\{ T_{t}\}$ is continuously $2$-differentiable in $t$, then $\ddot\lambda_j^\pm=\ddot \lambda_j^\pm(t)$ are unilaterally continuous, so that it holds that 
\[ \lim_{h\rightarrow \pm 0}\ddot \lambda_j^\pm(t+h)=\ddot \lambda_j^\pm(t) \] 
for any $t\in I$ and $j=1,2,\cdots$. 
\end{theorem}

\begin{theorem}\label{c2theorem} 
If $\{ T_t\}$ is continuously $2$-differentiable in $t$, the $C^1$ curves $\tilde C_j$, $k\leq j\leq k+m-1$ in Theorem \ref{rellich+} are $C^2$. 
\end{theorem}

Although the above theorems are to be extended to $p\geq 3$ of Definition \ref{def1}, we restrict ourselves to $p=1,2$ in this paper. Yet these results on $p=2$ are efficient to examine the harmonic concavity of $\lambda_j(t)$ in $t$ studied by \cite{gs52} for the first eigenvalue to (\ref{2}) with $\gamma_1=\emptyset$. We emphasize also that the treatment of the second derivatives is rather different from that of the first ones. See Remark \ref{rem9} in Section \ref{sec6}.

\section{Reduction to the abstract theory}\label{sec2}

 For the moment, we fix $t$ and treat the bi-Lipschitz homeomorphism $T=T_t:\Omega \rightarrow \Omega_t=T_t \Omega$. Let 
\[ \tilde \Omega=\Omega_t, \quad f(y)=g(x), \ y=Tx, \] 
and confirm the chain rule for this transformation of variables, that is,  
\[ \nabla g=\nabla f \ DT, \quad dy=(\det DT)dx \] 
for $\nabla g=(\frac{\partial g}{\partial x_1}, \cdots, \frac{\partial g}{\partial x_n})$, $\nabla f=(\frac{\partial f}{\partial y_1}, \cdots, \frac{\partial f}{\partial y_n})$, and 
\[ DT=\left( \begin{array}{ccc}
\frac{\partial y_1}{\partial x_1} & \cdots & \frac{\partial y_1}{\partial x_n} \\ 
\cdot & & \cdot \\ 
\frac{\partial y_n}{\partial x_1} & \cdots& \frac{\partial y_n}{\partial x_n} \end{array} \right). \] 

Putting $\tilde \gamma_i=\gamma_{it}=T\gamma_i$ for $i=0,1$, we take the eigenvalue problem  
\begin{equation} 
-\Delta u=\lambda u \ \mbox{in $\tilde \Omega$}, \quad u=0 \ \mbox{on $\tilde \gamma_0$}, \quad \frac{\partial u}{\partial \nu}=0 \ \mbox{on $\tilde\gamma_1$}, 
 \label{888}
\end{equation}  
that is, (\ref{7}) for $T=T_t$. Let 
\[ \tilde V=\{ v\in H^1(\tilde\Omega) \mid \left. v\right\vert_{\tilde \gamma_0}=0 \}, \] 
and introduce the weak form of (\ref{888}),  
\begin{equation} 
u\in \tilde V, \ \tilde A(u,v)=\lambda\tilde B(u,v), \quad \forall v\in \tilde V, 
 \label{10}
\end{equation} 
where  
\begin{equation} 
\tilde A(u,v)=\int_{\tilde\Omega} \nabla_y u\cdot \nabla_y v \ dy, \quad \tilde B(u,v)=\int_{\tilde\Omega} uv \ dy. 
 \label{12}
\end{equation}
 
Given $\phi\in V$, put 
\[ \psi(y)=\phi(x), \quad y=Tx. \] 
Then it holds that   
\[ \phi\in V \ \Leftrightarrow \ \psi\in \tilde V \] 
for $V\subset H^1(\Omega)$ defined by (\ref{v}). Writing 
\begin{equation} 
U(x)=u(y),  \  V(x)=v(y), \quad y=Tx, 
 \label{hoshi}
\end{equation} 
we obtain 
\begin{eqnarray*} 
\nabla_yu\cdot \nabla_yv & = & [\nabla_xU(DT)^{-1}]\cdot [\nabla_xV(DT)^{-1}] \\ 
& = & (\nabla_xU)(DT)^{-1}(DT)^{-1T}(\nabla_x V)^T \\ 
& = & (\nabla_xU)Q(\nabla_xV)^T=Q[\nabla_xU, \nabla_xV] 
\end{eqnarray*} 
for $Q=(DT)^{-1}(DT)^{-1T}$, where $F^T$ denotes the transpose of the matrix $F$. Then, (\ref{10}) means 
\[ \int_\Omega Q[\nabla_x U, \nabla_x V] a \ dx=\lambda\int_\Omega UVa \ dx \] 
for $a=\det DT$. The condition of normalization 
\[ \int_{\tilde \Omega}u^2 \ dy=1 \] 
is also transformed into  
\[ \int_\Omega U^2a \ dx=1. \] 

Under the family $\{T_t\}$ of homeomorphisms, therefore, the  weak form (\ref{88}) of (\ref{7}), is equivalent to (\ref{wf}) for $V\subset H^1(\Omega)$ defined by (\ref{v}) and $B_t$, $A_t$ defined by (\ref{ab})-(\ref{cd}). Here we confirm the following lemma. 

\begin{lemma}\label{equiv} 
The $j$-th eigenvalue of (\ref{7}) is equal to that of (\ref{wf}). 
\end{lemma} 
\begin{proof} 
For the moment, let $\tilde \lambda_j(t)$ be the $j$-th eigenvalues of (\ref{7}) and let $\lambda_j(t)$ be that of  (\ref{wf}) for (\ref{ab}) and (\ref{cd}). By the mini-max principle (\ref{minimax}), it holds that 
\begin{equation} 
\tilde \lambda_j(t)=\min_{\tilde L_j}\max_{v\in \tilde L_j\setminus \{0\}}\tilde R_t[v]=\max_{\tilde W_j}\min_{v\in \tilde W_j\setminus \{0\}}\tilde R_t[v], 
 \label{mm1} 
\end{equation} 
where $V_t\subset H^1(\Omega_t)$ is defined by (\ref{9}), 
\[ \tilde R_t[v]=\frac{\tilde A_t(v,v)}{\tilde B_t(v,v)}, \quad  
\tilde A_t(u,v)=\int_{\Omega_t}\nabla u\cdot \nabla v \ dx, \quad \tilde B_t(u,v)=\int_{\Omega_t}uv \ dx, \] 
and $\{ \tilde L_j\}$ and $\{ \tilde W_j\}$ are the families of all subspaces of $V_t$ with dimension and codimension $j$ and $j-1$, respectively.

 It holds also that 
\begin{equation} 
\lambda_j(t)=\min_{L_j}\max_{v\in L_j\setminus \{0\}} R_t[v]=\max_{W_j}\min_{v\in W_j\setminus \{0\}} R_t[v], 
 \label{mm2}
\end{equation} 
for  
\begin{equation} 
R_t[v]=\frac{A_t(v,v)}{B_t(v,v)},  
 \label{rt}
\end{equation} 
 and $\{ L_j\}$ and $\{ W_j\}$ denote the families of all subspaces of $V$ with dimension and codimension $j$ and $j-1$, respectively.

If the set $L$ is a $j$-dimensional subspace of $V$ there is $\phi_\ell \in L$, $1\leq \ell \leq j$, such that 
\[ \int_\Omega \phi_\ell \phi_{\ell'}dx=\delta_{\ell \ell'} \] 
and 
\[ \phi=\sum_{\ell=1}^jc_\ell \phi_\ell, \quad c_\ell=\int_\Omega \phi\phi_\ell \ dx  \] 
for any $\phi \in L$, which implies   
\[ \psi=\sum_{\ell=1}^jc_\ell \psi_\ell \] 
for $\psi=\phi\circ T_t^{-1}$ and $\psi_\ell=\phi\circ T_t^{-1}$. Hence we obtain $\mbox{dim} \ \tilde L_t\leq \mbox{dim} \ L$ for  
\[ \tilde L_t=\{ \phi\circ T_t^{-1} \mid \phi \in L\}. \] 

The reverse inequality follows similarly, and hence it holds that  
\[ \mbox{dim} \ \tilde L_t=\mbox{dim} \ L=j. \] 
Since $T_t:\Omega \rightarrow \Omega_t$ is a bi-Lipschitz homeomorphism, furthermore, $L\subset V$ if and only if $\tilde L_t\subset  V_t$. We thus obtain  
\begin{equation} 
\tilde \lambda_j(t)=\lambda_j(t)      
 \label{17}
\end{equation} 
by (\ref{mm1})-(\ref{mm2}).
\end{proof} 

We are ready to develop an abstract theory, writing $L^2$ norm in $X=L^2(\Omega)$ as $\vert \ \cdot \ \vert_X$. With $V$ in (\ref{v}), we recall that $\Vert \ \cdot \ \Vert_V$ denotes the norm in $V$ and that the inclusion $V\hookrightarrow X$ is compact. It holds also that 
\begin{equation} 
\vert v\vert_X\leq K\Vert v\Vert_V, \quad v\in V 
 \label{contie}
\end{equation} 
for $K=1$.

Henceforth, $C$ denotes the generic positive constant. 
The above $A_t:V\times V\rightarrow \R$ and $B_t:X\times X\rightarrow \R$ for $t\in I$ are symmetric bilinear forms, satisfying 
\begin{equation}  
\vert A_t(u,v)\vert \leq C\Vert u\Vert_V\Vert v\Vert_V, \ A_t(u,u)\geq \delta \Vert u\Vert_V^2, \quad u,v\in V 
 \label{coercive1}
\end{equation} 
and 
\begin{equation} 
\vert B_t(u,v)\vert \leq C\vert u\vert_X\vert v\vert_X, \ B_t(u,u)\geq \delta\vert u\vert_X^2, \quad u, v\in X 
 \label{coercive2}
\end{equation} 
for some $\delta>0$. Then the eigenvalues of (\ref{wf}) are denoted by 
\[ 0<\lambda_1(t)\leq \lambda_2(t)\leq \cdots \rightarrow +\infty. \] 
The weak and the strong convergences of $\{ u_j\}\subset Y$ to $u\in Y$ for $Y=X$ or $Y=V$ are, furthermore, indicated by 
\[ \mbox{w}\mathchar`-\lim_{j\rightarrow \infty}u_j=u  \ \mbox{in $Y$} \] 
and 
\[ \mbox{s}\mathchar`-\lim_{j\rightarrow \infty}u_j=u \ \mbox{in $Y$}, \] 
respectively.

\section{Continuity of eigenvalues and eigenspaces}\label{sec3}

Let $t\in I$ be fixed. We begin with the following theorem valid under the abstract setting in the previous section. 

\begin{theorem}\label{thm2}
The conditions   
\begin{eqnarray} 
& & \lim_{h\rightarrow 0}\sup_{\Vert u\Vert_V, \Vert v\Vert_V\leq 1}
\vert A_{t+h}(u,v)-A_t(u,v)\vert =0 \nonumber\\ 
& & \lim_{h\rightarrow 0}\sup_{\vert u\vert_X, \vert v\vert_X\leq 1}
\vert B_{t+h}(u,v)-B_t(u,v)\vert=0,  
 \label{18ab} 
\end{eqnarray}  
imply  
\begin{equation} 
\lim_{h\rightarrow 0}\lambda_j(t+h)=\lambda_j(t) 
 \label{16}
\end{equation} 
for any $j=1,2, \cdots$.    
\end{theorem} 

\begin{proof}
We note that the $j$-th eigenvalue of (\ref{7}) is given by the mini-max principle as in (\ref{mm2}), for the Rayleigh quotient $R_t[v]$ defined by (\ref{rt}).

Given $t$, let  
\begin{eqnarray} 
& & \alpha(h)=\sup_{\Vert u\Vert_V, \Vert v\Vert_V\leq 1}\vert (A_{t+h}-A_{t})(u,v)\vert \nonumber\\ 
& & \beta(h)=\sup_{\vert u\vert_X, \vert v\vert_X\leq 1}\vert (B_{t+h}-B_{t})(u,v)\vert. 
 \label{alphabeta}
\end{eqnarray}
We obtain  
\[ (A_{t+h}-A_{t})(v, v)\geq -\alpha(h)\Vert v\Vert_V^2\geq -\frac{\alpha(h)}{\delta}A_{t}(v,v) \] 
and 
\[ (A_{t+h}-A_{t})(v,v)\leq \alpha(h)\Vert v\Vert_V^2\leq \frac{\alpha(h)}{\delta}A_{t}(v,v) \] 
by (\ref{coercive1}) and (\ref{18ab}), 
which implies    
\[ (1-\delta^{-1}\alpha(h))A_{t}(v,v)\leq A_{t+h}(v,v)\leq (1+\delta^{-1}\alpha(h))A_{t}(v,v). \] 
Similarly, there arises that 
\[ (1-\delta^{-1}\beta(h))B_{t}(v,v)\leq B_{t+h}(v,v)\leq (1+\delta^{-1}\beta(h))B_{t}(v,v) \] 
for any $v\in V$. 

Then it follows that   
\[ (1-o(1))R_{t}[v]\leq R_{t+h}[v]\leq (1+o(1))R_{t}[v] \] 
uniformly in $v\in V\setminus \{0\}$, and hence 
\[ (1-o(1))\lambda_j(t)\leq \lambda_j(t+h)\leq (1+o(1))\lambda_j(t) \] 
by (\ref{18ab}). Thus we obtain (\ref{16}). 
\end{proof}

Let $u_j(t)\in V$ be the eigenfunction of  (\ref{wf}) corresponding to the eigenvalue $\lambda_j(t)$:  
\begin{equation} 
B_t(u_j(t), u_{j'}(t))=\delta_{jj'}, \quad A_t(u_j(t), v)=\lambda_j(t)B_t(u_j(t),v), \ \forall v\in V.
 \label{19}
\end{equation} 
Fix $t$, assume (\ref{18}), and define $\lambda$ by (\ref{put}). Although this multiplicity $m$ is not stable under the perturbation of $t$, we obtain the following theorem concerning the continuity of eigenspaces with respect to $t$.  

Let 
\begin{equation} 
Y_t^\lambda=\langle u_j(t) \mid k\leq j\leq k+m-1 \rangle 
 \label{ym}
\end{equation} 
be the subspace of $X$ generated by the above eigenfunctions $u_j(t)$ for $k\leq j\leq k+m-1$.

\begin{lemma}\label{thm3}
Under the above situation, any $h_\ell\rightarrow 0$ admits a subsequence, denoted by the same symbol, such that the limits  
\begin{equation} 
\mbox{s}\mathchar`-\lim_{\ell\rightarrow \infty}u_j(t+h_\ell)=\phi_j\in Y_t^\lambda \ \mbox{in $V$}, \quad  k\leq j\leq k+m-1 
 \label{strong}
\end{equation} 
exist. In particular it holds that  
\begin{equation}  
B_{t}(\phi_j, \phi_{j'})=\delta_{jj'}, \quad k\leq j, j'\leq k+m-1. 
 \label{bb}
\end{equation} 
\end{lemma} 

\begin{proof} 
We have 
\begin{equation} 
\lim_{h\rightarrow 0}\lambda_j(t+h)=\lambda_j(t), \quad k\leq j\leq k+m-1 
 \label{20}
\end{equation} 
by Theorem \ref{thm2}. It holds also that 
\begin{equation} 
\Vert u_j(t')\Vert_V\leq C, \quad k\leq j\leq k+m-1, \ \vert t'\vert<\varepsilon_0  
 \label{23}
\end{equation} 
by (\ref{contie})-(\ref{coercive2}) and (\ref{19}): 
\[ \delta \Vert u_j(t')\Vert_V^2\leq A_{t'}(u_j(t'), u_j(t'))=\lambda_j(t')B_t(u_j(t'), u_j(t'))=\lambda_j(t'). \] 

Given $h_\ell\rightarrow 0$, therefore, we have a subsequence, denoted by the same symbol, which admits the weak limits   
\begin{equation} 
\mbox{w}\mathchar`-\lim_{\ell\rightarrow \infty}u_j(t+h_\ell)=\phi_j \ \mbox{in $V$}, \quad k\leq j\leq k+m-1, 
 \label{wc}
\end{equation}  
for some $\phi_j\in V$. From (\ref{19}) it follows that 
\[ A_{t}(\phi_j, v)=\lambda_j(t)B_{t}(\phi_j, v), \quad \forall v\in V \] 
and hence 
\[ \phi_j\in Y_t^\lambda, \quad k\leq j\leq k+m-1. \] 

Since $V\hookrightarrow X$ is compact, the weak convergence (\ref{wc}) implies the strong convergence
\begin{equation}  
\mbox{s}\mathchar`-\lim_{\ell\rightarrow \infty}u_j(t+h_\ell)=\phi_j \ \mbox{in $X$},  
 \label{l2strong}
\end{equation} 
and hence (\ref{bb}). Now we improve this weak convergence (\ref{wc}) to the strong convergence (\ref{strong}) in $V$, using (\ref{18ab}). 

For this purpose, we put  
\[ v=u_j(t+h_\ell)-\phi_j \] 
and recall (\ref{put}) and (\ref{20}). Then there arises that  
\begin{eqnarray*} 
\delta \Vert v\Vert_V^2 & \leq & A_{t+h_\ell}(u_j(t+h_\ell)-\phi_j, v) \nonumber\\ 
& = & A_{t+h_\ell}(u_j(t+h_\ell), v) - (A_{t+h_\ell}-A_{t})(\phi_j, v)-A_{t}(\phi_j, v) \nonumber\\ 
& = & \lambda_j(t+h_\ell)B_{t+h_\ell}(u_j(t+h_\ell), v) - (A_{t+h_\ell}-A_{t})(\phi_j,v) -\lambda_j(t)B_{t}(\phi_j, v) \nonumber\\ 
& = & (\lambda_j(t+h_\ell)-\lambda_j(t))B_{t+h_\ell}(u_j(t+h_\ell),v) \nonumber\\ 
& & +\lambda_j(t)(B_{t+h_\ell}-B_{t})(u_j(t+h_\ell),v) +\lambda_j(t)B_{t}(u_j(t+h_\ell)-\phi_j,v) \nonumber\\ 
& & -(A_{t+h_\ell}-A_{t})(\phi_j, v) \nonumber\\  
& \leq &  C\vert \lambda_j(t+h_\ell)-\lambda_j(t)\vert \cdot \vert u_j(t+h_\ell)\vert_X\cdot K\Vert v\Vert_V \nonumber\\ 
& & + \lambda_j(t)\cdot \beta(h_\ell)\cdot \vert u_j(t+h_\ell)\vert_X\cdot K\Vert v\Vert_V \nonumber\\ 
& & +\lambda_j(t)\cdot C\vert u_j(t+h_\ell)-\phi_j\vert_X\cdot K\Vert v\Vert_V +\alpha(h_\ell)\Vert \phi_j\Vert_V\cdot \Vert v\Vert_V, 
\end{eqnarray*} 
and hence 
\begin{eqnarray*} 
\Vert v \Vert_V & \leq & C\{ \vert \lambda_j(t+h_\ell)-\lambda_j(t)\vert +\beta(h_\ell)+\vert u_j(t+h_\ell)-\phi_j\vert_X +\alpha(h_\ell)\} \\ 
& =& o(1) 
\end{eqnarray*}
by (\ref{18ab}) and (\ref{l2strong}).  
\end{proof} 

\begin{remark}\label{ambiguity} 
If $m=1$, the eigenfunction $u_j(t)$ in (\ref{26}) is uniquely determined by (\ref{19}) up to the multiplication of $\pm 1$, which implies $\phi_j=\pm u_j(t)$. In the other case of $m\geq 2$, the eigenfunction which attains (\ref{17}) does not satisfy this property. Hence we have   
\[ \phi_j=\sum_{i=1}^mq^i_ju_i(t) \] 
for $Q=(q^i_j)$ satisfying $Q^{T}=Q^{-1}$ in Theorem \ref{thm3}. In other words, the eigenfunction $u_j(t)$ corresponding to $\lambda_j(t)$ has more varieties than $\pm 1$ multiplication, although the eigenspace $Y_t^\lambda$ is determined. By this ambiguity the limits $\phi_j$ in (\ref{strong}) depend on the sequence $h_\ell\rightarrow 0$, which makes the argument below to be complicated. 
\end{remark} 

\section{First derivatives}\label{sec5}

If $\{ T_t\}$ is $1$-differentiable in the setting of Section \ref{sec1}, we can put   
\[ \dot A_t(u,v)=\int_\Omega \dot Q_t[\nabla u, \nabla v]a_t+Q_t[\nabla u, \nabla v]\dot a_t \ dx, \quad u, v\in V \] 
and 
\[ \dot B_t(u,v)=\int_\Omega uv \dot a_t \ dx, \quad u, v\in X. \] 
These $\dot A_t:V\times V\rightarrow \R$ and $\dot B_t:X\times X\rightarrow \R$ in (\ref{ab}) are bilinear forms  
satisfying 
\begin{eqnarray} 
& & \vert \dot A_t(u,v)\vert \leq C\Vert u\Vert_V\Vert v\Vert_V, \quad u,v\in V \nonumber\\ 
& & \vert \dot B_t(u,v)\vert \leq C\vert u\vert_X\vert v\vert_X, \quad u, v\in X 
 \label{d11}  
\end{eqnarray} 
and   
\begin{eqnarray} 
& &  \lim_{h\rightarrow 0}\frac{1}{h}\sup_{\Vert u\Vert_V, \Vert v\Vert_V\leq 1}
 \left\vert \left(A_{t+h}-A_t-h\dot A_t\right)(u,v)\right\vert =0 \nonumber\\ 
& & \lim_{h\rightarrow 0}\frac{1}{h}\sup_{\vert u\vert_X, \vert v\vert_X\leq 1}
\left\vert \left(B_{t+h}-B_t-h\dot B_t\right)(u,v)\right\vert =0.
 \label{abd}
\end{eqnarray} 
Hence Theorem \ref{thm2.1} in Section \ref{sec2} is reduced to the following abstract theorem. 
In this theorem, the assumption made by Theorem \ref{thm2} is valid, and therefore, there arises that (\ref{16}). 

\begin{theorem}\label{thm2.1abst} 
Let $X, V$ be Hilbert spaces over $\mathbb{R}$, with compact embedding $V\hookrightarrow X$. Let  $A_t:V\times V\rightarrow \mathbb{R}$ and $B_t:X\times X\rightarrow \mathbb{R}$ be symmetric bilinear forms satisfying (\ref{coercive1})-(\ref{coercive2}) for any $t\in I$. Given $t$, assume the existence of the bilinear forms $\dot A_t:V\times V\rightarrow \mathbb{R}$ and $\dot B_t:X\times X\rightarrow \mathbb{R}$ such that (\ref{d11})-(\ref{abd}). Assume, finally, (\ref{18}) with $\ell=k$ and $k=m$ for $\lambda_j(t)$, $k\leq k+m-1$, defined by (\ref{mm2})-(\ref{rt}). Then the conclusion of Theorem \ref{thm2.1} holds.  
\end{theorem}

Before proceeding to the rigorous proof, we develop a formal argument, writing (\ref{wf}) as 
\begin{equation} 
u_t\in V, \ B_t(u_t,u_t)=1, \quad A_t(u_t,v)=\lambda_tB_t(u_t,v), \ v\in V.  
 \label{b1}
\end{equation} 
First, taking a formal differentiation in $t$ in this equality, we obtain   
\begin{equation} 
\dot A_t(u_t, v)+A_t(\dot u_t, v)=\dot \lambda_tB_t(u_t, v)+\lambda_t\dot B_t(u_t, v) + \lambda_tB_t(\dot u_t, v), \quad \forall v\in V. 
 \label{25} 
\end{equation} 
Putting $v=u_t$, we obtain 
\begin{equation} 
\dot A_t(u_t,u_t)+A_t(\dot u_t, u_t)=\dot\lambda_t+\lambda_t\dot B_t(u_t, u_t)+\lambda_tB_t(\dot u_t, u_t)   
 \label{26}
\end{equation} 
by 
\begin{equation} 
B_t(u_t,u_t)=1. 
 \label{27}
\end{equation} 

To eliminate $\dot u_t$ in (\ref{26}), second,  we use (\ref{27}) to deduce 
\begin{equation} 
\dot B_t(u_t,u_t)+2B_t(\dot u_t, u_t)=0. 
 \label{28}
\end{equation} 
From 
\[ \lambda_t=\lambda_tB_t(u_t,u_t)=A_t(u_t, u_t) \] 
it is derived also that 
\begin{equation} 
\dot \lambda_t=\dot A_t(u_t,u_t)+2A_t(\dot u_t, u_t) 
 \label{29}
\end{equation} 
and then (\ref{26}) is replaced by 
\[ \dot A_t(u_t,u_t)+\frac{1}{2}\dot\lambda_t-\frac{1}{2}\dot A_t(u_t, u_t)=\dot \lambda_t+\lambda_t\dot B_t(u_t, u_t)-\frac{1}{2}\lambda_t\dot B_t(u_t,u_t), \]  
or, 
\begin{equation} 
\dot \lambda_t=\dot A_t(u_t, u_t)-\lambda_t\dot B_t(u_t, u_t). 
 \label{30}
\end{equation} 
As is noticed in Remark \ref{ambiguity} in Section \ref{sec3}, if the eigenspace 
\[ Y_t^\lambda=\langle u_j(t) \mid k\leq j\leq k+m-1\rangle \] 
corresponding to the eigenvalue $\lambda=\lambda_j(t)$, $k\leq j\leq k+m-1$, to (\ref{wf}),  is one-dimensional as in $m=1$, the eigenfunction $u_t$ in (\ref{b1}) is unique up to the multiplication of $\pm 1$, and this ambiguity is canceled in (\ref{30}). This property is valid even if $m\geq 2$ as in Remark \ref{orthgonalt}. 

Turning to the rigorous proof, we use the following lemma, recalling $u_j(t)\in V$, $Y^\lambda_t$, and $\phi_j$ in (\ref{19}), (\ref{ym}), and (\ref{strong}), respectively. 

\begin{lemma}\label{lem4}
Under the assumption of Theorem \ref{thm2.1abst}, any $h_\ell\rightarrow 0$ admits a subsequence, denoted by the same symbol, such that  
\begin{equation} 
\lim_{\ell\rightarrow\infty}\frac{1}{h_\ell}\{ \lambda_j(t+h_\ell)-\lambda_j(t)\} =\dot A_{t}(\phi_j, \phi_j)-\lambda_j(t)\dot B_{t}(\phi_j, \phi_j)
 \label{wakewake}
\end{equation} 
for $k\leq j\leq k+m-1$. It holds also that 
\begin{equation} 
\dot A_{t}(\phi_j, \phi_{j'})-\lambda_j(t)\dot B_{t}(\phi_j, \phi_{j'})=0, \quad k\leq j\neq j'\leq k+m-1. 
 \label{jj}
\end{equation} 
\end{lemma} 

\begin{proof} 
Let $k\leq j, j'\leq k+m-1$. Since 
\begin{eqnarray*} 
& & A_{t+h}(u_j(t+h)-\phi_j, u_{j'}(t+h)-\phi_{j'}) = A_{t+h}(u_j(t+h), u_{j'}(t+h)) \\ 
& & \quad -A_{t+h}(u_j(t+h), \phi_{j'}) -A_{t+h}(\phi_j, u_{j'}(t+h))+A_{t+h}(\phi_j, \phi_{j'}) 
\end{eqnarray*} 
and 
\begin{eqnarray*} 
& & A_{t}(u_j(t+h)-\phi_j, u_{j'}(t+h)-\phi_{j'})=A_{t}(u_j(t+h), u_{j'}(t+h)) \\ 
& & \quad -A_{t}(u_j(t+h), \phi_{j'})-A_{t}(\phi_j, u_{j'}(t+h))+A_{t}(\phi_j,\phi_{j'}), 
\end{eqnarray*} 
it holds that  
\begin{eqnarray} 
& & h(A_{t+h}-A_{t})\left(\frac{u_j(t+h)-\phi_j}{h}, \frac{u_{j'}(t+h)-\phi_{j'}}{h}\right) \nonumber\\
& & \quad =\frac{1}{h}(A_{t+h}-A_{t})(u_j(t+h)-\phi_j, u_{j'}(t+h)-\phi_{j'}) \nonumber\\ 
& & \quad = \frac{1}{h}(A_{t+h}-A_{t})(u_j(t+h), u_{j'}(t+h)) 
+\frac{1}{h}(A_{t+h}-A_{t})(\phi_j, \phi_{j'}) \nonumber\\ 
& & \qquad -\frac{1}{h}(A_{t+h}-A_t)(u_j(t+h), \phi_{j'})-\frac{1}{h}(A_{t+h}-A_t)(\phi_j, u_{j'}(t+h)).  
 \label{33}
\end{eqnarray} 

In (\ref{33}), first, we obtain  
\[ 
\frac{1}{h}(A_{t+h}-A_{t})(u_j(t+h), u_{j'}(t+h))=\dot A_{t}(u_j(t+h),u_{j'}(t+h))+o(1) \] 
as $h\rightarrow 0$ by (\ref{23}) and (\ref{abd}). Then (\ref{strong}) implies 
\begin{eqnarray} 
& & \lim_{\ell\rightarrow\infty}\frac{1}{h_\ell}(A_{t+h_\ell}-A_{t})(u_j(t+h_\ell), u_{j'}(t+h_\ell)) \nonumber\\ 
& & \quad =\lim_{\ell\rightarrow \infty}\dot A_t(u_j(t+h_\ell), u_{j'}(t+h_\ell))
=\dot A_t(\phi_j, \phi_{j'}) 
 \label{35} 
\end{eqnarray} 
and also 
\begin{equation}  
\lim_{\ell\rightarrow\infty}\frac{1}{h_\ell}(A_{t+h_\ell}-A_{t})(\phi_j, \phi_{j'})
= \dot A_{t}(\phi_j, \phi_{j'}). 
 \label{35+}
\end{equation}

Second, it holds that 
\begin{eqnarray*} 
& & A_{t+h}(u_j(t+h), \phi_{j'})=\lambda_j(t+h)B_{t+h}(u_j(t+h), \phi_{j'}) \\ 
& & A_{t}(u_j(t+h), \phi_{j'})=\lambda_j(t)B_{t}(u_j(t+h), \phi_{j'}) 
\end{eqnarray*} 
by $\phi_{j'}\in {Y_t^\lambda}$, which implies     
\begin{eqnarray} 
& & \frac{1}{h_\ell}(A_{t+h_\ell}-A_{t})(u_j(t+h_\ell), \phi_{j'}) \nonumber\\ 
& & \quad =\frac{1}{h_\ell}(\lambda_j(t+h_\ell)B_{t+h_\ell}(u_j(t+h_\ell), \phi_{j'})-\lambda_j(t)B_{t}(u_j(t+h_\ell), \phi_{j'})) \nonumber\\ 
& & \quad =\frac{1}{h_\ell}(\lambda_j(t+h_\ell)-\lambda_j(t))B_{t+h_\ell}(u_j(t+h_\ell), \phi_{j'}) \nonumber\\ 
& & \qquad +\frac{1}{h_\ell}\lambda_j(t)(B_{t+h_\ell}-B_{t})(u_j(t+h_\ell),\phi_{j'}) \nonumber\\ 
& & \quad =  \frac{1}{h_\ell}(\lambda_j(t+h_\ell)-\lambda_j(t))(\delta_{jj'}+o(1))+ \lambda_j(t)\dot B_{t}(\phi_j, \phi_{j'})+o(1) 
\end{eqnarray} 
by 
\begin{eqnarray*} 
B_{t+h_\ell}(u_j(t+h_\ell), \phi_{j'}) & = & B_{t}(u_j(t+h_\ell), \phi_{j'})+o(1) \\ 
& = & B_t(\phi_j, \phi_{j'})+o(1) =\delta_{jj'}+o(1) 
\end{eqnarray*} 
and (\ref{abd}). Similarly, it follows that   
\begin{eqnarray} 
& & \frac{1}{h_\ell}(A_{t+h_\ell}-A_{t})(\phi_j, u_{j'}(t+h_\ell)) \nonumber\\ 
& & \quad =\frac{1}{h_\ell}(\lambda_j(t+h_\ell)-\lambda_j(t))(\delta_{jj'}+o(1))+
\lambda_j(t)\dot B_{t}(\phi_j, \phi_{j'})+o(1).  
 \label{43}
\end{eqnarray} 

Finally, we obtain   
\begin{eqnarray} 
& & \left\vert h_\ell(A_{t+h_\ell}-A_{t})\left(\frac{u_j(t+h_\ell)-\phi_j}{h_\ell}, \frac{u_{j'}(t+h_\ell)-\phi_{j'}}{h_\ell}\right)\right\vert \nonumber\\ 
& & \quad \leq Ch_\ell^2\left\Vert \frac{u_j(t+h_\ell)-\phi_j}{h_\ell}\right\Vert_V\cdot \left\Vert \frac{u_{j'}(t+h_\ell)-\phi_{j'}}{h_\ell}\right\Vert_V \nonumber\\ 
& & \quad = C\Vert u_j(t+h_\ell)-\phi_j\Vert_V\cdot \Vert u_{j'}(t+h_\ell)-\phi_{j'}\Vert_V =o(1) 
 \label{45x}
\end{eqnarray} 
by (\ref{strong}) and (\ref{abd}). Then equalities (\ref{wakewake})-(\ref{jj}) follow from (\ref{35})-(\ref{45x}) as 
\begin{eqnarray*} 
0 & = & 2\dot A_t(\phi_j, \phi_{j'})-\frac{2}{h_\ell}(\lambda_j(t+h_\ell)-\lambda_j(t))(\delta_{jj'}+o(1)) \\ 
& & -2\lambda_j(t)\dot B_t(\phi_j, \phi_{j'})+o(1), \quad \ell\rightarrow \infty.  
\end{eqnarray*} .    
\end{proof} 

Below we confirm that the process of taking subsequence in the previous lemma is not necessary, if $h_\ell\rightarrow 0$ is unilateral as in $h_\ell\rightarrow +0$ or $h_\ell\rightarrow -0$. Theorem \ref{thm2.1} is thus reduced to the following theorem.  

\begin{theorem}\label{thm4}
Under the assumption of Theorem \ref{thm2.1abst}, the unilateral limits    
\begin{equation} 
\dot \lambda_j^\pm(t)=\lim_{h\rightarrow \pm 0}\frac{1}{h}\{ \lambda_j(t+h)-\lambda_j(t)\} 
 \label{49} 
\end{equation} 
exist, and it holds that 
\begin{equation}  
\dot \lambda^+_j(t)=\mu_{j-k+1}, \ \dot\lambda_j^-(t)=\mu_{k+m-j}, \quad k\leq j\leq k+m-1. 
 \label{pm}
\end{equation} 
Here, $\mu_{q}$, $1\leq q\leq m$, is the $q$-th eigenvalue of 
\begin{equation} 
u\in Y_t^\lambda, \quad E^\lambda_{t}(u,v)=\mu B_t(u,v), \ \forall v\in Y_t^\lambda,  
 \label{51}
\end{equation} 
where $Y_t^\lambda$ is the $m$-dimensional eigenspace of (\ref{wf}) corresponding to the eigenvalue $\lambda$ of  (\ref{put}) defined by (\ref{ym}), and   
\begin{equation} 
E^\lambda_{t}=\dot A_{t}-\lambda \dot B_{t}. 
 \label{et}
\end{equation} 
In particular, it holds that 
\[ \dot \lambda_j^+(t)=\dot \lambda_{2k+m-1-j}^-(t), \quad k\leq j\leq k+m-1. \] 
\end{theorem} 

\begin{proof} 
Since $Y_t^\lambda$ is $m$-dimensional, the eigenvalue problem (\ref{51}) admits $m$-eigenvalues denoted by   
\[ \mu_1\leq \cdots \leq \mu_m. \] 
By Lemma \ref{lem4}, on the other hand, any $h_\ell\rightarrow 0$ takes a subsequence, denoted by the same symbol, satisfying (\ref{wakewake})-(\ref{jj}) for some 
\[ \phi_j\in Y_t^\lambda, \quad k\leq j\leq k+m-1, \] 
with (\ref{bb}). 

This lemma ensures also the existence of  
\begin{equation} 
\tilde \mu_j=\lim_{\ell\rightarrow\infty}\frac{1}{h_\ell}(\lambda_j(t+h_\ell)-\lambda_j(t)),  
 \label{limit}
\end{equation} 
and the equalities 
\begin{equation} 
E^\lambda_t(\phi_j, \phi_{j'})=\delta_{jj'}\tilde \mu_j, \quad k\leq j, j' \leq k+m-1. 
 \label{68}
\end{equation} 
We thus obtain    
\[ \phi_j\in {Y_t^\lambda}, \quad B_{t}(\phi_j, \phi_j)=1, \quad E^\lambda_{t}(\phi_j,v)=\tilde \mu_jB_{t}(\phi_j, v), \ \forall v\in Y_t^\lambda, \] 
and therefore, $\mu=\tilde \mu_j$ is an eigenvalue of (\ref{51}). 

If $h_\ell\rightarrow +0$, there arises that  
\[ \tilde \mu_k\leq \cdots \leq \tilde \mu_{k+m-1} \] 
by  
\[ \lambda_{k}(t+h)\leq \cdots \leq \lambda_{k+m-1}(t+h),  \] 
and hence 
\[ \tilde \mu_j=\mu_{j-k+1}, \quad k\leq j\leq k+m-1. \] 
Then we obtain the result because the value $\tilde \mu_j$ in (\ref{limit}) is independent of the sequence $h_\ell\rightarrow +0$. 

In the other case of $h_\ell\rightarrow -0$, we obtain 
\[ \tilde\mu_j=\mu_{k+m-j}, \quad k\leq j\leq k+m-1,  \] 
and the result follows similarly. 
\end{proof} 

Theorem \ref{thm2+} is reduced to the following abstract theorem. 

\begin{theorem}\label{abstconti}
Let the assumption of Theorem \ref{thm2.1abst} hold in $I$. Fix $t\in I$, and assume  
\begin{eqnarray} 
& & \lim_{h\rightarrow 0}\sup_{\Vert u\Vert_V, \Vert v\Vert_V\leq 1}\left\vert \dot A_{t+h}(u,v)-\dot A_t(u,v)\right\vert=0 \nonumber\\ 
& & \lim_{h\rightarrow 0}\sup_{\vert u\vert_X, \vert v\vert_X\leq 1}\left \vert \dot B_{t+h}(u,v)-\dot B_t(u,v)\right\vert=0.   
 \label{72x}
\end{eqnarray} 
Then, it follows that  
\[ \lim_{h\rightarrow \pm 0}\dot \lambda_j^\pm(t+h)=\dot \lambda_j^\pm (t). \] 
\end{theorem} 

\begin{proof} 
Assume (\ref{18})-(\ref{put}) and take $k\leq j\leq k+m-1$. Since the assumption of Theorem \ref{thm2.1abst} holds in $I$, any $t'\in I$ admits $u_j(t')\in V$ such that  
\begin{equation} 
B_{t'}(u_j(t'), u_j(t'))=1, \ A_{t'}(u_j(t'), u_j(t'))=\lambda_j(t')B_{t'}(u_j(t'), u_j(t'))  
 \label{73x}
\end{equation} 
and 
\begin{equation} 
\dot \lambda_j^+(t')=\dot A_{t'}(u_j(t'), u_j(t'))-\lambda_j(t')\dot B_{t'}(u_j(t'), u_j(t')).  
 \label{73xx}
\end{equation} 
by Lemma \ref{lem4} and Theorem \ref{thm4}.

Given $t$ in this theorem, and take $h_\ell\rightarrow +0$ and $u_j(t')$ in (\ref{73x}) for $t'=t+h_\ell$.  Hence  there is a subsequence denoted by the same symbol such that (\ref{strong}) with $\phi_j\in V$. Then it holds that    
\[ \dot\lambda_j^+(t)=\dot A_t(\phi_j, \phi_j)-\lambda_j(t)\dot B_t(\phi_j, \phi_j). \] 
We thus obtain 
\begin{eqnarray*} 
\lim_{\ell\rightarrow \infty}
\dot \lambda_j^+(t+h_\ell) & = & \lim_{\ell\rightarrow \infty}
\{ \dot A_{t+h_\ell}(u_j(t+h_\ell), u_j(t+h_\ell)) \\ 
& & -\lambda_j(t+h_\ell)\dot B_{t+h_\ell}(u_j(t+h_\ell), u_j(t+h_\ell)) \} \\ 
& = & \dot  A_t(\phi_j, \phi_j)-\lambda_j(t)\dot B_t(\phi_j, \phi_j)=\dot \lambda_j^+(t) 
\end{eqnarray*} 
by (\ref{73xx}), and hence 
\[ \lim_{h\rightarrow +0}\dot \lambda_j^+(t+h)=\dot \lambda_j^+(t) \] 
because $h_\ell\rightarrow +0$ is arbitrary. 

The proof of 
\[ \lim_{h\rightarrow -0}\dot \lambda_j^-(t+h)=\dot \lambda_j^-(t) \] 
is similar. 
\end{proof}

\begin{remark}\label{uniform1} 
The limits (\ref{limit}) exist for any $j$ under the conditions
(\ref{coercive1}), (\ref{coercive2}), (\ref{d11}), and  (\ref{abd}). If
these conditions are satisfied for any $t\in (-\varepsilon_0,
\varepsilon_0)$, the limits (\ref{49}) are unilaterally
 locally uniform in $t\in I$.
In fact, if not, there are, for example, $t_k\downarrow  t_0\in$
$(-\varepsilon_0, \varepsilon_0)$,
 $\delta>0$, and $h_\ell\rightarrow 0$, such that 
\[ \left\vert \frac{1}{h_\ell}(\lambda_j(t_k+h_\ell)-\lambda_j(t_k))-\dot \lambda_j^+(t_k)\right\vert \geq \delta. \] 
Then we obtain 
\[ \left\vert
 \frac{1}{h_\ell}(\lambda_j(t_0+h_\ell)-\lambda_j(t_0))
 -\dot \lambda_j^+(t_0)\right\vert \geq \delta, \]
a contradiction with $\ell\rightarrow \infty$. 
\end{remark}

\section{Rearrangement of eigenvalues}\label{sec:rearrange}

For simplicity we introduce the following notations to prove Theorem \ref{rellich+}. Recall $I=(-\varepsilon_0, \varepsilon_0)$, and let $f_j\in C^0(I)$, $1\leq j\leq m$, satisfy  
\begin{equation} 
f_1(t)\leq \cdots \leq f_m(t), \quad t\in I. 
 \label{p0}
\end{equation} 
Assume the existence of the unilateral limits 
\begin{equation} 
\dot f_j^\pm(t)=\lim_{h\rightarrow \pm 0}\frac{1}{h}(f_j(t+h)-f_j(t)) 
 \label{p1}
\end{equation} 
and 
\begin{equation} 
\lim_{h\rightarrow \pm 0}\dot f_j(t+h)=\dot f_j^\pm(t) 
 \label{p2}
\end{equation} 
for any $j$ and $t$. Assume, finally, 
\begin{equation} 
\dot f_j^+(t)=\dot f_{2k+n-j-1}^-(t), \quad k\leq j\leq k+n-1,  
 \label{p3}
\end{equation} 
provided that 
\begin{equation} 
f_{k-1}(t)<f_k(t)=\cdots = f_{k+n-1}(t)<f_{k+n}(t),  
 \label{cluster}
\end{equation}
where $1\leq n\leq m$, $1\leq k\leq m-n+1$, and $t\in I$. 
In (\ref{cluster}) we understand 
\[ f_0(t)=-\infty, \quad f_{m+1}(t)=+\infty. \]  

We call 
\[ K=K_{k,n}(t)=\{ k, \cdots, k+n-1\} \] 
the {\it $n$-cluster} at $t$ with entry $k$ if (\ref{cluster}) arises, and also   
\[ p(K)=\max\{ j \mid k\leq j\leq k+[\frac{n}{2}]-1, \ \dot f_j^+(t)<\dot f_{2k+n-1-j}^+(t)\} -k+1  \] 
its {\it $p$-value}. Here, we understand $p(K)=0$ if 
\[ \dot f_k^+(t)=\dot f_{k+n-1}^+(t), \] 
noting 
\[ k\leq i\leq j\leq k+n-1 \quad \Rightarrow \quad \dot f_i^+(t)\leq \dot f_j^+(t). \]   

We construct a rearrangement of $C^0$-curves
\[ C_j=\{ f_j(t) \mid t\in I\}, \quad 1\leq j\leq m,  \] 
denoted by $\tilde C_j$, $1\leq j\leq m$, so that are $C^1$ in $t\in I$.  This rearrangement is done only on 
\[ I_1=\{ t\in I \mid \mbox{there exists a cluster $K$ at $t$ such that $p(K)\geq 1$}.\}. \] 

To introduce this rearrangement, we note the following facts in advance. First, given $2\leq n\leq m$, let    
\[ I_1^n=\{ t\in I \mid \mbox{there exists an $n$-cluster $K$ at $t$ such that $p(K)\geq 1$}\}.  \] 
If $t\in I_1^n$ and $K=K_{k,n}(t)$ satisfies $p(K)\geq 1$, it holds that 
\begin{equation} 
f_{k+n-1}(t')>f_k(t'), \quad 0<\vert t'-t\vert \ll 1.  
 \label{bunretsu}
\end{equation} 
Hence this $t$ is an isolated point of $I_1^n$. In particular, each $I_1^n$, $2\leq n\leq m$, is at most countable, and hence so is $I_1$ by   
\[ I_1=\bigcup_{n=2}^mI_1^n. \] 
Second, given $t\in I_1$ and $1\leq j\leq m$, if $j\in K=K_{k,n}(t)$ holds for some $K$ in $p(K)\geq 1$, this $K$ is unique.  

\begin{definition}\label{transversal}  
The curves $\tilde C_j$, $1\leq j\leq m$, are called the transversal rearrangement of $C_j $, $1\leq j\leq m$, if the following operations are done. 
\begin{enumerate} 
\item If $t\in I_1$, $1\leq j\leq m$, and $j\in K$ hold for $K=K_{k,n}(t)$ with $p(K)\geq 1$, the curve $C_j$ for $k\leq j\leq k+p-1$ and $k-n-p\leq j\leq k-n-1$ on the right, is connected to $C_{2k+n-j-1}$ on the left at $t$, where $p=p(K)$. 
\item No rearrangements to $C_j$, $1\leq j\leq m$, are done otherwise.  
\end{enumerate} 
\end{definition} 

The curves $\tilde C_j$, $1\leq j\leq m$, are uniquely constructed from $C_j$, $1\leq j\leq m$, by this transversal rearrangement. From  the results in the previous section, Theorem \ref{rellich+} is a consequence of the following theorem. 

\begin{theorem}\label{c1theorem}
Under the above situation, the $C^0$ curves $\tilde C_j$, $1\leq j\leq m$, made by the transversal rearrangement of $C_j$, $1\leq j\leq m$, are $C^1$.  
\end{theorem} 

\begin{proof} 
This theorem is obvious if $m=1$. Now we show it by an induction on $m$, assuming the assertion up to $m-1$. 

Take $t_0\in I\setminus \overline{I_1}$, and make the transversal rearrangement of $C_j$, $1\leq j\leq q$, toward left and right diections.  Let $t_\ell$, $\ell=1,2,\cdots$, be the successive points of $I_1$ in the left direction: 
\[ t_\ell \in I_1, \ t_{\ell-1}>t_\ell, \ (t_\ell, t_{\ell-1})\cap I_1=\emptyset, \quad \ell=1,2,\cdots. \] 
If $\{ t_\ell \}$ is finite, these $C_j$'s are successfully rearranged to $C^1$ curves on $(-\varepsilon_0, t_0)$. If not, there is 
\[ t_\ast=\lim_{\ell\rightarrow \infty}t_\ell \in [-\varepsilon_0, t_0]. \] 
Then the case $t_\ast=-\varepsilon_0$ ensures the same conclusion. 

Letting $t_\ast>-\varepsilon_0$, we show that $\{ \tilde C_j\mid 1\leq j\leq m\}$ are $C^1$ curves on $(t_\ast-\delta, t_\ast+\delta)$ for $0<\delta \ll 1$.  Once this fact is proven, we can repeat this process up to $t=-\varepsilon_0$ by a covering argument. Turning to the right direction, we conclude that $\tilde C_j$, $1\leq j\leq m$, are $C^1$ curves on $I=(-\varepsilon_0, \varepsilon_0)$.  

To this end we distinguish the cases $t_\ast\in I_1$ and $t_\ast \not\in I_1$. 

If $t_\ast \in I_1$, first, the above  assertion follows from the assumption of induction. In fact, each $K=K_{k,n}(t_\ast)$ with $p(K)\geq 1$ admits (\ref{bunretsu}), while 
\[ \tilde f_j(t)=\left\{ \begin{array}{ll} 
f_j(t), & 0<t-t_\ast\ll 1\\ 
f_{2k+n-j-1}(t), & 0<t_\ast-t\ll 1\end{array} \right.,  \quad j=k, \ j=k+n-1 \] 
are $C^1$ around $t=t_\ast$. Then we apply the assumption of induction to $C_j$, $k+1\leq j\leq k+n-2$, to get $n$-$C^1$ curves in $(t_\ast-\delta, t_\ast+\delta)$ made by $C_j$, $k\leq j\leq k+n-1$. Operating this process to any $K=K_{k,n}(t_\ast)$ with $p(K)\geq 1$ at $t=t_\ast$, we get $C^1$ curves $\tilde C_j$, $1\leq j\leq m$, in $(t_\ast-\delta, t_\ast+\delta)$ by this transversal rearrangement of $C_j$, $1\leq j\leq m$, at $t=t_\ast$.  

If $t_\ast\not \in I_1$, second, we take the cluster decomposition of $\{ 1,\cdots, m\}$ at $t=t_\ast$, that is,  
\[ 1=k_1<k_1+n_1=k_2<\cdots< k_{s-1}+n_{s-1}=k_s<k_s+n_s=m \] 
satisfying  
\[ \bigcup_{r=1}^sK_{k_r, n_r}(t_\ast)=\{ 1, \cdots, m\}. \] 
Since 
\[ p(K_{k_r, n_r}(t_\ast))=0, \quad 1\leq r\leq s \] 
holds by the assumption, there are $a_r$, $1\leq r\leq s$, such that 
\[ \dot f_j^+(t_\ast)=a_r, \quad \forall j\in K_{k_r, n_r}. \]

We obtain, on the other hand,  
\[ f_{k_{r-1}}(t)<f_{k_r}(t), \ \vert t-t_\ast\vert\ll 1, \ 1\leq r\leq s+1 \] 
under the agreement 
\[ f_{k_0}(t)=-\infty, \quad f_{k_{s+1}}(t)=+\infty, \] 
and hence $\tilde C_j$, $1\leq j\leq m$, are $C^1$ on $[t_\ast, t_\ast+\delta)$ for $0<\delta \ll 1$. If $t_\ast$ is not a right accumulating point of $I_1$, therefore, these $\tilde C_j$, $1\leq j\leq m$, made by the transversal rearrangement of $C_j$, $1\leq j\leq m$, are $C^1$ on $(t_\ast-\delta, t_\ast+\delta)$. 

In the other case that $t_\ast$ is a right accumulating point of $I_1$, these $\tilde C_j$, $1\leq j\leq m$, are $C^1$ on $(t_\ast-\delta, t_\ast]$, similarly. Then it holds that  
\[ \dot f_j^-(t_\ast)=a_r, \quad \forall j\in K_{k_r, n_r}=K_{k_r,n_r}(t_\ast) \] 
by $p(K_{k_r,n_r})=0$, and hence $\tilde C_j$, $1\leq j\leq m$, are $C^1$ on $(t_\ast-\delta, t_\ast+\delta)$. 
\end{proof}

\section{Second derivatives} \label{sec6}

If $T_t:\Omega\rightarrow \Omega_t$ is $2$-differentiable, we have the other bilinear forms $\ddot A_t:V\times V\rightarrow \R$ and $\ddot B_t:X\times X\rightarrow \R$ satisfying 
\begin{eqnarray} 
& & \vert \ddot A_t(u,v)\vert \leq C\Vert u\Vert_V\Vert v\Vert_V, \quad u,v\in V \nonumber\\ 
& & \vert \ddot B_t(u,v)\vert \leq C\vert u\vert_X\vert v\vert_X, \quad \ u,v\in X 
 \label{abd2+} 
\end{eqnarray} 
uniformly in $t$ and 
\begin{eqnarray} 
& & \lim_{h\rightarrow 0}\frac{1}{h^2}\sup_{\Vert u\Vert_V, \Vert v\Vert_V\leq 1}\left\vert \left(A_{t+h}-A_t-h\dot A_t-\frac{h^2}{2}\ddot A_t\right)(u,v)\right\vert=0 \nonumber\\ 
& & \lim_{h\rightarrow 0} \frac{1}{h^2}\sup_{\vert u\vert_X, \vert v\vert_X \leq 1}\left\vert \left( B_{t+h}-B_t-h\dot B_t-\frac{h^2}{2}\ddot B_t\right)(u,v)\right\vert =0 
 \label{abd2}
\end{eqnarray} 
for each $t$. Hence Theorem \ref{thmsecond} is reduced to the following abstract theorem. 

\begin{theorem}\label{secondabst}
In Theorem \ref{thm2.1abst}, assume, futhermore, (\ref{abd2+})-(\ref{abd2}). Then the conclusion of Theorem \ref{thmsecond} holds. 
\end{theorem}

For the moment we develop a formal argument as in the first derivative. Assuming (\ref{30}), first, we deduce  
\begin{eqnarray} 
\ddot\lambda_t & = & \ddot A_t(u_t, u_t)+2\dot A_t(\dot u_t, u_t)-\dot \lambda_t\dot B_t(u_t,u_t)-\lambda_t\ddot B_t(u_t, u_t)-2\lambda_t\dot B_t(\dot u_t, u_t) \nonumber\\ 
& = & 2(\dot A_t-\lambda_t\dot B_t)(\dot u_t, u_t)+D_t(u_t, u_t) 
 \label{39} 
\end{eqnarray} 
for  
\begin{equation} 
D_t(u,v)=\ddot A_t(u,v)-\dot \lambda_t\dot B_t(u,v)-\lambda_t\ddot B_t(u,v), \quad u, v\in V. 
 \label{dt}
\end{equation} 
Putting $v=\dot u_t$ in (\ref{25}), second, we reach 
\begin{equation} 
(\dot A_t-\lambda_t\dot B_t)(u_t, \dot u_t)=-(A_t-\lambda_tB_t)(\dot u_t, \dot u_t)+\dot\lambda_tB_t(u_t, \dot u_t). 
 \label{40}
\end{equation} 
Then, (\ref{28}), (\ref{39}), and (\ref{40}) imply 
\begin{eqnarray} 
& & \ddot \lambda_t = -2(A_t-\lambda_tB_t)(z_\ast, z_\ast)+2\dot \lambda_tB_t(u_t,z_\ast)+D_t(u_t, u_t) \nonumber\\ 
& & = -2(A_t-\lambda_tB_t)(z_\ast, z_\ast)-\dot\lambda_t\dot B_t(u_t, u_t)+D_t(u_t, u_t) \nonumber\\ 
& & = -2(A_t-\lambda_tB_t)(z_\ast, z_\ast)+\ddot A_t(u_t, u_t)-2\dot \lambda_t\dot B_t(u_t,u_t)-\lambda_t\ddot B_t(u_t, u_t)  
 \label{41} 
\end{eqnarray} 
for $z_\ast=\dot u_t$. 

We observe that $\dot u_t\in V$ is not uniquely determined by (\ref{25}), which is derive formally also. It has, more precisely, the ambiguity of addition of an element in $Y_t^\lambda$. This ambiguity, however, cancels in  (\ref{41}) by equality  (\ref{vanish2}) below. 

We now develop a rigorous argument valid even to (\ref{18}). Define $\lambda$ by (\ref{put}) and recall $Y_t^\lambda=\langle u_j(t) \mid k\leq j\leq k+m-1\rangle$ for $u_j(t)$ satisfying (\ref{19}). Let, also, 
\[ C_{t'}^j=A_{t'}-\lambda_j(t')B_{t'}, \quad k\leq j\leq k+m-1, \ t'\in I. \] 
Then we obtain 
\[ C_t^j=A_t-\lambda B_t\equiv C_t \] 
and    
\begin{equation} 
C_t(u,v)=0, \quad \forall (u,v) \in Y_t^\lambda\times V.    
 \label{vanish2}
\end{equation} 

By Lemma \ref{lem4}, given $h_\ell\rightarrow 0$, which may not be of definite sign, we have a subsequence, denoted by the same symbol, satisfying  (\ref{strong}), 
\begin{equation}  
\mbox{s}\mathchar`-\lim_{\ell\rightarrow\infty}u_j(t+h_\ell)=\phi_j\in
 Y_t^\lambda  \ \mbox{in $V$}, \quad k\leq j\leq k+m-1.   
 \label{phi1}
\end{equation}  
There exists  
\begin{equation} 
\dot\lambda_j^\ast=\lim_{\ell\rightarrow\infty}\frac{1}{h_\ell}(\lambda_j(t+h_\ell)-\lambda_j(t)), \quad k\leq j\leq k+m-1 
 \label{67}
\end{equation} 
passing to a subsequence, with the equality  
\[ \dot\lambda_j^\ast\delta_{jj'}=\dot A_t(\phi_j, \phi_{j'})-\lambda\dot B_t(\phi_j, \phi_{j'}), \quad k\leq j, j'\leq k+m-1.  \] 

\begin{remark}\label{remast}
If we take the subsequence of $h_\ell\rightarrow 0$ to be unilateral, the limit $\dot \lambda_j^\ast$ in (\ref{67}) exists and is either $\dot \lambda_j^+(t)$ or $\dot\lambda_j^-(t)$.  
\end{remark} 

Let 
\[ \dot C_{t}^{\ast j}=\dot A_t-\dot\lambda_j^\ast B_t-\lambda \dot B_t.  \] 
It holds that  
\[ \lim_{\ell\rightarrow \infty}\sup_{\Vert u\Vert_V \leq 1, \Vert v\Vert_V \leq 1}\left\vert \frac{1}{h_\ell}(C_{t+h_\ell}^j-C_t)(u,v)-\dot C_{t}^{\ast j}(u,v)\right\vert =0, \quad k\leq j\leq k+m-1 \] 
by (\ref{abd}), and also   
\begin{equation} 
\dot C_{t}^{\ast j}(u,v)=0, \quad \forall u, v\in Y_t^\lambda, \ k\leq j\leq k+m-1   
 \label{vanish}
\end{equation} 
by Lemma \ref{lem4}. Let 
\[  
z_{\ell}^j=\frac{1}{h_\ell}(u_j(t+h_\ell)-\phi_j). 
\]   

\begin{lemma}\label{lem6}
It holds that    
\begin{equation} 
\lim_{\ell\rightarrow \infty}C_t(z_{\ell}^j,v)=-\dot C_{t}^{\ast j}(\phi_j, v), \quad \forall v\in V.   
 \label{ct}
\end{equation} 
\end{lemma}
\begin{proof} 
Given $v\in V$, we obtain  
\begin{eqnarray*} 
C_t(z_{\ell}^j,v) & = & \frac{1}{h_\ell}C_t(u_j(t+h_\ell)-\phi_j, v) = \frac{1}{h_\ell}C_t(u_j(t+h_\ell), v) \\ 
& = & -\frac{1}{h_\ell}(C_{t+h_\ell}^j-C_t)(u_j(t+h_\ell),v)=-\dot C_{t}^{\ast j}(u_j(t+h_\ell), v)+o(1) 
\end{eqnarray*} 
by (\ref{23}) and (\ref{vanish2}). Then (\ref{ct}) follows from (\ref{strong}). 
\end{proof} 

Recall that 
\[ R:X\rightarrow Y_t^\lambda=\langle u_j(t) \mid k\leq j\leq k+m-1\rangle \] 
is the orthogonal projection with respect to $B_t(\cdot, \cdot)$, and $P=I-R$. There is a unique $z_\ast^j \in PV$ satisfying 
\begin{equation}  
C_t(z_\ast^j, v)=-\dot C_{t}^{\ast j}(\phi_j, v), \quad \forall v\in PV. 
 \label{70e}
\end{equation} 

\begin{remark}\label{rempp} 
Equality (\ref{70e}) ensures 
\begin{equation} 
\gamma_{\dot \lambda_j^\ast}(\phi_j)=z_\ast^j, \quad k\leq j\leq k+m-1 
 \label{phi2}
\end{equation} 
under the notation of Definition \ref{def2}. 
\end{remark} 

\begin{lemma}\label{lem7}  
It holds that 
\begin{equation} 
\mbox{w}\mathchar`-\lim_{\ell\rightarrow \infty}Pz_{\ell}^j=z_\ast^j \quad \mbox{in $V$}, \quad  
k\leq j\leq k+m-1. 
 \label{pw}
\end{equation} 
\end{lemma} 
\begin{proof} 
Lemma \ref{lem6} ensures that $\{ Pz_{\ell}^j\}$ converges weakly in $PV$, and hence is bounded there:  
\[ \Vert Pz_{\ell}^j\Vert_V\leq C. \] 
Then, passing to a subsequence denoted by the same symbol, there is $\tilde z_j\in PV$ such that 
\[ \mbox{w}\mathchar`-\lim_{\ell\rightarrow \infty}Pz_{\ell}^j=\tilde z_j, \] 
which satisfies 
\begin{equation} 
C_t(\tilde z_j, v)=-\dot C_{t}^{\ast j}(\phi_j, v), \quad \forall v\in V 
 \label{72e}
\end{equation} 
by Lemma \ref{lem6}. Since such $\tilde z_j\in PV$ is unique, we obtain the result with $\tilde z_j=z_\ast^j$.  
\end{proof} 

\begin{remark}\label{remplus}
Since (\ref{72e}) holds with $\tilde z_j=z^j_\ast$, this $z_j^\ast \in PV$ defined by (\ref{70e}) satisfies 
\[ 
C_t(z_\ast^j, v)=-\dot C_{t}^{\ast j}(\phi_j, v), \quad \forall v\in V. 
\] 

\end{remark} 

\begin{remark}\label{rem9} 
Generally, the inequality 
\[ \Vert z_\ell^j\Vert_V\leq C \] 
is not expected to hold, which causes the other difficulty in later arguments. In fact, if $m=1$ and $\phi_j=u_j(t)$, for example, this property means  
\[ \vert B_t(u_j(t), w_{h_\ell}^j)\vert \leq C, \quad w_h^j=\frac{1}{h}(u_j(t+h)-u_j(t)). \] 
In the formal argument, we have actually (\ref{28}), which, however, does not assure the actual convergence   
\begin{equation} 
\lim_{h\rightarrow 0}B_t(u_j(t), w_h^j)=-\frac{1}{2}\dot B_t(u_j(t), u_j(t)).  
 \label{50}
\end{equation} 
In fact, the equality  
\[ 1=B_{t+h}(u_j(t+h), u_j(t+h))=B_t(u_j(t), u_j(t)) \] 
just implies  
\begin{eqnarray*} 
0 & = & \frac{1}{h}\{ B_{t+h}(u_j(t+h), u_j(t+h))-B_t(u_j(t),u_j(t))\} \\ 
& = & \frac{1}{h}(B_{t+h}-B_t)(u_j(t+h), u_j(t+h)) \\ 
& & +\frac{1}{h}\{ B_t(u_j(t+h), u_j(t+h))-B_t(u_j(t), u_j(t))\} \\ 
& = & \dot B_t(u_j(t),u_j(t))+o(1)+\frac{1}{h}B_t(u_j(t+h)+u_j(t), u_j(t+h)-u_j(t)) \\ 
& = & \dot B_t(u_j(t),u_j(t))+2B_t(\frac{u_j(t+h)+u_j(t)}{2}, w_h^j) +o(1)
\end{eqnarray*} 
and hence 
\begin{equation} 
\lim_{h\rightarrow 0}B_t(\frac{u_j(t+h)+u_j(t)}{2}, w_h^j)=-\frac{1}{2}\dot B_t(u_j(t), u_j(t)) 
 \label{left}
\end{equation}   
differently from (\ref{50}).  Here, the condition $\Vert w_h^j\Vert_X=O(1)$ is necessary to conclude 
\[ B_t(u_j(t+h), w_h^j)=B_t(u_j(t), w_h^j)+o(1) \] 
in the left-hand side of (\ref{left}) from $\phi_j=u_j(t)$ in (\ref{strong}). Our purpose, however, was to assure $\Vert w_h^j\Vert_V=O(1)$, which is reduced to $\Vert w_h^j\Vert_X=O(1)$ by $\Vert P_hw_h^j\Vert_V=O(1)$. This is a circular reasoning.   
\end{remark}

\begin{lemma}\label{lem9}
It holds that 
\begin{equation} 
\mbox{s}\mathchar`-\lim_{\ell\rightarrow \infty}Pz_{\ell}^j=z_\ast^j \quad \mbox{in $V$}, \qquad k\leq j\leq k+m-1. 
 \label{phi3}
\end{equation} 
\end{lemma} 

\begin{proof} 
Since $V\hookrightarrow X$ is compact, we have  
\begin{equation} 
\mbox{s}\mathchar`-\lim_{\ell\rightarrow\infty} Pz_{\ell}^j=z_\ast^j \quad \mbox{in $X$}  
 \label{53}
\end{equation} 
in the previous lemma. Then we obtain 
\begin{eqnarray*} 
\delta \Vert Pz_{\ell}^j-z_\ast^j\Vert_V^2 & \leq & A_t(Pz_{\ell}^j-z_\ast^j, Pz_{\ell}^j-z_\ast^j)= C_t(Pz_{\ell}^j-z_\ast^j, Pz_{\ell}^j-z_\ast^j)+o(1) \\ 
& = & C_t(Pz_{\ell}^j, Pz_{\ell}^j-z_\ast^j)+o(1) =C_t(z_{\ell}^j, Pz_{\ell}^j-z_\ast^j)+o(1)  
\end{eqnarray*} 
by (\ref{pw}) and (\ref{53}). 

Since $\phi_j\in Y_t^\lambda$ it holds that  
\begin{eqnarray*} 
C_t(z_{\ell}^j, Pz_{\ell}^j-z_\ast^j) & = & \frac{1}{h_\ell}C_t(u_j(t+h_\ell), Pz_{\ell}^j-z_\ast^j) \\
& = & \frac{1}{h_\ell}(C_t-C_{t+h_\ell}^j)(u_j(t+h_\ell), Pz_{\ell}^j-z_\ast^j) \\ 
& = & -\dot C_{t}^{\ast j}(u_j(t+h_\ell), Pz_{\ell}^j-z_\ast^j)+o(1)  \\ 
& = &  -\dot C_{t}^{\ast j}(\phi_j, Pz_{\ell}^j-z_\ast^j)+o(1)=o(1) 
\end{eqnarray*} 
by (\ref{strong}) and (\ref{pw}), because 
\[ v\in V \ \mapsto \ \dot C_{t}^{\ast j}(\phi_j, v)\in \R \] 
is a bounded linear mapping. Then the result follows as 
\[ \lim_{\ell\rightarrow \infty}\Vert Pz_\ell^j-z_\ast^j\Vert_V=0. \]  
\end{proof} 

\begin{remark} 
The limit $z_\ast^j$ in (\ref{phi3}) depends on the sequence $h_\ell\rightarrow 0$ because it is prescribed by (\ref{phi1}) and (\ref{phi2}). The limit $\ddot \lambda_j^\ast$ in the following lemma depends also $h_\ell\rightarrow 0$, but this ambiguity is canceled unilaterally. Hence these limits are uniquely determined as  either $h\downarrow 0$ or $h\uparrow 0$, because they are characterized by an eigenvalue problem on a finite dimsnsional space, similarly to the first derivative of $\lambda_j(t)$.  See Theorem \ref{thm14} below.  
\end{remark}

\begin{lemma}\label{lem12} 
There exists 
\begin{equation} 
\ddot \lambda_j^\ast\equiv \lim_{\ell\rightarrow \infty}\frac{2}{h_\ell^2}(\lambda_j(t+h_\ell)-\lambda_j(t)-h_\ell\dot\lambda_j^\ast), \quad 
k\leq j\leq k+1-m 
 \label{h2}
\end{equation} 
with 
\[ \ddot \lambda_j^\ast =(\ddot A_t-\lambda \ddot B_t-2\dot \lambda_j^\ast\dot B_t)(\phi_j, \phi_j)-2C_t(z_\ast^j, z_\ast^j).  \] 
\end{lemma} 

\begin{proof} 
By (\ref{vanish}) and Lemma \ref{lem9}, we have  
\begin{eqnarray} 
C_t(z_\ell^j, z_\ell^j) & = & C_t(Pz_\ell^j, Pz_\ell^j)=C_t(z_\ast^j, z_\ast^j)+o(1) \nonumber\\ 
& = & C_t(Pz_\ell^j, z_\ast^j)+o(1)=C_t(z_\ell^j, z_\ast^j)+o(1). 
 \label{kome0}
\end{eqnarray} 
It holds that   
\begin{eqnarray*} 
C_t(z_\ell^j, z_\ast^j) & = & \frac{1}{h_\ell}C_t(u_j(t+h_\ell)-\phi_j, z_\ast^j)=\frac{1}{h_\ell}C_t(u_j(t+h_\ell), z_\ast^j) \\ 
& = & \frac{1}{h_\ell}(C_t-C_{t+h_\ell}^j)(u_j(t+h_\ell), z_\ast^j) = -\dot C_{t}^{\ast j}(u_j(t+h_\ell), z_\ast^j)+o(1) \\ 
& = & -\dot C_{t}^{\ast j}(\phi_j, z_\ast^j)+o(1)  
\end{eqnarray*} 
by (\ref{vanish2}) and $\phi_j\in Y_t^\lambda$, which implies 
\begin{equation} 
C_t(z_\ast^j, z_\ast^j)=-\dot C_t^{\ast j}(\phi_j, z_\ast^j)+o(1)
 \label{kome}
\end{equation} 
by (\ref{kome0}). It holds also that 
\begin{eqnarray} 
\dot C_{t}^{\ast j}(\phi_j, z_\ast^j) & = & \frac{1}{h_\ell}\dot C_{t}^{\ast j}(\phi_j, P(u_j(t+h_\ell)-\phi_j))+o(1) \nonumber\\ 
& = & \frac{1}{h_\ell}\dot C_{t}^{\ast j}(\phi_j, u_j(t+h_\ell)-\phi_j)+o(1) \nonumber\\ 
& = & \frac{1}{h_\ell}\dot C_{t}^{\ast j}(\phi_j, u_j(t+h_\ell))+o(1) 
 \label{kome2} 
\end{eqnarray} 
by (\ref{vanish}) and $\phi_j\in Y_t^\lambda$.   

Here, we use the asymptotics  
\begin{eqnarray*} 
& & C_{t+h_\ell}^j(\phi_j, u_j(t+h_\ell)) = C_{t}^j(\phi_j, u_j(t+h_\ell))+h_\ell\dot C_{t}^{\ast j}(\phi_j, u_j(t+h_\ell)) \\ 
& & +\frac{1}{2}h_\ell^2\ddot C_{t,\ell}^{\ast j}(\phi_j, u_j(t+h_\ell))+o(h_\ell^2) 
\end{eqnarray*} 
for  
\begin{equation} 
\ddot C_{t,\ell}^{\ast j}=\ddot A_t-\lambda \ddot B_t-2\dot \lambda_j^\ast\dot B_t-\frac{2}{h_\ell^2}(\lambda_j(t+h_\ell)-\lambda_j(t)-h\dot \lambda_j^\ast)B_t, 
 \label{kome3}
\end{equation} 
derived from (\ref{abd}) and (\ref{abd2}). Since 
\[ C_{t+h_\ell}^j(\phi_j, u_j(t+h_\ell))=C_t^j(\phi_j, u_j(t+h_\ell))=0 \] 
holds by (\ref{vanish2}), we obtain 
\begin{eqnarray*} 
C_t^j(z_\ast^j, z_\ast^j) & = & - \dot C_{t}^{\ast j}(\phi_j, z_\ast^j)+o(1)=-\frac{1}{h_\ell}\dot C_{t}^{\ast j}(\phi_j, u_j(t+h_\ell))+o(1) \\ 
& = & \frac{1}{2}\ddot C_{t,\ell}^{\ast j}(\phi_j, u_j(t+h_\ell))+o(1) = \frac{1}{2}\ddot C_{t,\ell}^{\ast j}(\phi_j, \phi_j)+o(1) \\ 
& = & \frac{1}{2}(\ddot A_t-\lambda \ddot B_t-2\dot \lambda_j^\ast\dot B_t)(\phi_j, \phi_j) -\frac{1}{h_\ell^2}(\lambda_j(t+h_\ell)-\lambda_j(t)-h_\ell\dot \lambda_j^\ast)+o(1) 
\end{eqnarray*} 
by (\ref{kome})-(\ref{kome2}) and $B_t(\phi_j, \phi_j)=1$. Then it follows that   
\[ 
\lim_{\ell\rightarrow \infty}\frac{2}{h_\ell^2}(\lambda_j(t+h_\ell)-\lambda_j(t)-h_\ell\dot \lambda_j^\ast) \\ 
=(\ddot A_t-\lambda\ddot B_t-2\dot\lambda_j^\ast\dot B_t)(\phi_j, \phi_j)-2C_t(z_\ast^j, z_\ast^j),   
\]  
and the proof is complete. 
\end{proof} 

\begin{lemma}\label{lem13}
If 
$\dot\lambda_\ast\equiv \dot \lambda_j^\ast=\dot \lambda_{j'}^\ast$ arises for some $k\leq j\neq j'\leq k+m-1$, 
then it holds that 
\[ (\ddot A_t-\lambda\ddot B_t-2\dot\lambda_\ast \dot B_t)(\phi_j, \phi_{j'})=2C_t(z_\ast^j, z_\ast^{j'}). \] 
\end{lemma} 

\begin{proof} 
As in the previous lemma we obtain    
\begin{eqnarray*} 
 C_t(z_\ast^j, z_\ast^{j'}) & = & C_t(Pz_\ell^j, z_\ast^{j'})+o(1)=C_t(z_\ell^j, z_\ast^{j'})+o(1) \\ 
& = & \frac{1}{h_\ell}C_t(u_j(t+h_\ell)-\phi_j, z_\ast^{j'})+o(1)=\frac{1}{h_\ell}C_t(u_j(t+h_\ell), z_\ast^{j'})+o(1) \\ 
& = & \frac{1}{h_\ell}(C_t-C_{t+h_\ell})(u_j(t+h), z_\ast^{j'})+o(1) \\ 
& = & -\dot C_{t}^{\ast j}(u_j(t+h_\ell), z_\ast^{j'})+o(1)=-\dot C_{t\ast}^j(\phi_j, z_\ast^{j'})+o(1) \\ 
& = & -\frac{1}{h_\ell}\dot C_{t}^{\ast j}(\phi_j, P(u_{j'}(t+h_\ell)-\phi_{j'}))+o(1) \\ 
& = & -\frac{1}{h_\ell}\dot C_{t}^{\ast j}(\phi_j, u_{j'}(t+h_\ell)-\phi_{j'})+o(1) \\ 
& = & -\frac{1}{h_\ell}\dot C_{t}^{\ast j}(\phi_j, u_{j'}(t+h_\ell))+o(1)  
\end{eqnarray*} 
by $\phi_j, \phi_{j'}\in Y_t^\lambda$. Then it holds that 
\begin{eqnarray*} 
& & C_{t+h_\ell}^j(\phi_j, u_{j'}(t+h_\ell))=C_t(\phi_j, u_{j'}(t+h_\ell))+h_\ell \dot C_{t}^{\ast j}(\phi_j, u_{j'}(t+h_\ell)) \\ 
& & \quad +\frac{1}{2}h_\ell^2 \ddot C_{t,\ell}^{\ast j}(\phi_j, u_{j'}(t+h_\ell))+o(h_\ell^2) 
\end{eqnarray*} 
with 
\[ C_t(\phi_j, u_{j'}(t+h_\ell))=0 \] 
and 
\begin{eqnarray*} 
C^j_{t+h_\ell}(\phi_j, u_{j'}(t+h_\ell)) & = & (C_{t+h_\ell}^j-C_{t+h_\ell}^{j'})(\phi_j, u_{j'}(t+h_\ell)) \\ 
& = & (\lambda_{j'}(t+h_\ell)-\lambda_j(t+h_\ell))B_{t+h_\ell}(\phi_j, u_{j'}(t+h_\ell)) 
\end{eqnarray*} 
by (\ref{vanish2}), to conclude 
\begin{eqnarray} 
C_t(z_\ast^j, z_\ast^{j'}) & = & \frac{1}{h_\ell^2}(\lambda_{j}(t+h_\ell)-\lambda_{j'}(t+h_\ell))B_{t+h_\ell}(\phi_j, u_{j'}(t+h_\ell)) \nonumber\\ 
& & +\frac{1}{2}\ddot C_{t,\ell}^{\ast j}(\phi_j, u_{j'}(t+h_\ell))+o(1). 
 \label{csecond}  
\end{eqnarray} 

Then we use  
\begin{eqnarray*} 
& & \lambda_j(t+h_\ell)=\lambda+h_\ell \dot \lambda_j^\ast+\frac{h_\ell^2}{2}\ddot \lambda_j^\ast+o(h_\ell^2) \\ 
& & \lambda_{j'}(t+h_\ell)=\lambda+h_\ell \dot \lambda_{j'}^\ast+\frac{h_\ell^2}{2}\ddot \lambda_{j'}^\ast+o(h_\ell^2) 
\end{eqnarray*}
with $\dot \lambda_j^\ast=\dot\lambda_{j'}^\ast$, to deduce  
\begin{eqnarray*} 
& & \lim_{\ell\rightarrow
 \infty}\frac{1}{h_\ell^2}(\lambda_j(t+h_\ell)-\lambda_{j'}(t+h_\ell))B_{t+h_\ell}(\phi_j,  u_{j'}(t+h_\ell)) \\  
& & \quad = \frac{1}{2}(\ddot \lambda_j^\ast+\ddot \lambda_{j'}^\ast)B_t(\phi_j, \phi_{j'})=0.  
\end{eqnarray*} 
Then the result follows from (\ref{kome3}), (\ref{csecond}), and the previous lemma. 
\end{proof} 

Recall $F^{\lambda, \lambda'}_t$ in Definition \ref{def2}. 

\begin{lemma} 
Define  
$\tilde \mu_k\leq \cdots \leq \tilde\mu_{k+m-1}$ 
by 
\[ \{ \tilde \mu_j \mid k\leq j\leq k+m-1\}=\{ \dot \lambda_j^\ast \mid k\leq j\leq k+m-1\}, \] 
and assume $k\leq \ell<r\leq k+m$ be such that  
\[ \tilde \mu_{\ell-1}<\tilde \mu\equiv \tilde \mu_{\ell}=\cdots = \tilde \mu_{r-1}<\tilde \mu_{r}. \] 
under the agreement of 
\[ \tilde \mu_{k-1}=-\infty, \quad \tilde\mu_{k+m}=+\infty. \] 
Then, $\sigma=\ddot \lambda_j^\ast$, $\ell\leq j\leq r-1$, is an eigenvalue of  
\[ u\in Y_{\lambda,t}^{\ell, r}, \quad F^{\lambda, \tilde \mu}_t(u,v)=\sigma B_t(u,v), \ \forall v\in Y_{\lambda,t}^{\ell, r}, \]
for 
\[ Y_{\lambda,t}^{\ell,r}=\langle u_j(t) \mid \ell\leq j\leq r-1\rangle \subset Y^\lambda_t. \] 
\end{lemma} 

\begin{proof} 
Lemmas \ref{lem12} and \ref{lem13} imply 
\[ F^{\lambda, \tilde \mu}_t(\phi_j, \phi_{j'})=\delta_{jj'}\ddot \lambda_j^\ast, \quad \ell\leq j, j'\leq r-1, \] 
and hence the result follows from  
\[ Y_{\lambda,t}^{\ell,r}=\langle \phi_j \mid \ell\leq j\leq r-1\rangle, \quad B_t(\phi_j, \phi_{j'})=\delta_{jj'}. \]  
\end{proof} 

Theorem \ref{secondabst} is now reduced to the following theorem. 

\begin{theorem}\label{thm14} 
Fix $t\in I$, and assume (\ref{18}) for $\ell=k$ and $n=m$. Put (\ref{put}) and let $k\leq \ell<r\leq k+m$ be such that 
\begin{equation} 
\dot\lambda^+_{\ell-1}(t)<\lambda'\equiv \dot \lambda^+_{\ell}(t)=\cdots=\dot \lambda_{r-1}^+(t)<\dot \lambda_{r}^+(t).   
 \label{equal}
\end{equation} 
Then there exists 
\begin{equation} 
\lambda_j''(t)=\lim_{h\rightarrow +0}\frac{2}{h^2}(\lambda_j(t+h)-\lambda-h\lambda'), \quad \ell \leq j\leq r-1.  
 \label{bilateral}
\end{equation} 
It holds, furthermore, that 
\begin{equation} 
\lambda_j''(t)=\sigma_{j-\ell+1}, \quad \ell\leq j\leq r-1,  
 \label{105}
\end{equation} 
where $\sigma_q$, $1\leq q\leq r-\ell$, denotes the $q$-th eigenvalue of  
\[ u\in Y_{\lambda,t}^{\ell, r}, \qquad F_t^{\lambda, \lambda'}(u,v)=\sigma B_t(u,v), \ \forall v\in Y_{\lambda,t}^{\ell, r}. \] 
 If 
\[ \dot\lambda^-_{\ell-1}(t)>\lambda'\equiv \dot \lambda^-_{\ell}(t)=\cdots=\dot \lambda_{r-1}^-(t)>\dot \lambda_{r}^-(t),  \] 
there arises that  
\[ \lambda_j''(t)=\lim_{h\rightarrow -0}\frac{2}{h^2}(\lambda_j(t+h)-\lambda-h\lambda'), \quad \ell \leq j\leq r-1 \] 
with (\ref{105}). 
\end{theorem} 
\begin{proof} 
In the previous lemma, we obtain 
\[ \ddot \lambda_\ell^\ast\leq \cdots \leq \ddot \lambda_{r-1}^\ast. \] 
Hence the result follows similarly to Theorem \ref{thm4}.
\end{proof} 

\begin{remark}\label{uniform2}
By the above theorem and Remark \ref{remast}, the limits (\ref{18+}),  
\[ \ddot\lambda_j^\pm(t)=\lim_{h\rightarrow \pm 0}\frac{2}{h^2}(\lambda_j(t+h)-\lambda_j(t)-h\dot \lambda_j^\pm(t)) \] 
exist for any $j$, provided that the conditions (\ref{coercive1}), (\ref{coercive2}), (\ref{d11}), (\ref{abd}), (\ref{abd2+}), and (\ref{abd2}) hold. If these conditions hold for any $t$ locally uniformly in $I$, these limits are locally uniform in $t\in I$. 
\end{remark}

Theorem \ref{secondconti} is reduced to the following abstract theorem. The proof is similar to that of Theorem \ref{abstconti}.  

\begin{theorem}\label{thm25}
Let the assumption of Theorem \ref{secondabst} hold for any $t$. Fix $t\in I$, and assume, furthermore, 
\begin{eqnarray} 
& & \lim_{h\rightarrow 0}\sup_{\Vert u\Vert_V, \Vert v\Vert_V\leq 1}\left\vert \ddot A_{t+h}(u,v)-\ddot A_t(u,v)\right\vert=0 \nonumber\\ 
& & \lim_{h\rightarrow 0}\sup_{\vert u\vert_X, \vert v\vert_X\leq 1}\left\vert \ddot B_{t+h}(u,v)-\ddot B_t(u,v)\right\vert=0. 
 \label{abdot}
\end{eqnarray} 
Then it holds that 
\[ \lim_{h\rightarrow \pm 0}\ddot \lambda_j^\pm(t+h)=\ddot \lambda_j^\pm(t). \] 
\end{theorem}

Finally, Theorem \ref{c2theorem} is reduced to the following abstract theorem. 

\begin{theorem}\label{finaltheorem} 
If (\ref{abdot}) is valid to any $t$ in the previous theorem, $C^1$ curves $\tilde C_j$, $1\leq j\leq m$, in Theorem \ref{rellich+} are $C^2$. 
\end{theorem} 
 
\begin{proof} 
Define $\tilde \lambda_j(t)$ by
\[ \tilde C_j=\{ \tilde \lambda_j(t) \mid t\in I\}, \ 1\leq j\leq m.  \]
From the proof of Theorem \ref{c1theorem}, Theorems \ref{thm14} and \ref{thm25} guarantee the existence of 
\begin{equation}   
\tilde \lambda_j''(t)=\lim_{h\rightarrow 0}\frac{1}{h^2}(\tilde \lambda_j(t+h)-\tilde \lambda_j(t)-h\tilde \lambda_j'(t)) 
 \label{last} 
\end{equation} 
together with its continuity in $t$, 
\[ \lim_{h\rightarrow \pm 0}\tilde \lambda_j''(t+h)=\tilde \lambda_j''(t) \] 
for any $t$ and $j$. This convergence (\ref{last}), furthermore, is locally uniform in $t\in I$ by Remark \ref{uniform2}. 

Then it follows that 
\begin{eqnarray*} 
& & \tilde \lambda_j(t+h)=\tilde \lambda_j(t)+h \tilde\lambda_j'(t)+\frac{h^2}{2}\tilde \lambda_j''(t)+o(h^2)  \\ 
& & \tilde \lambda_j(t)=\tilde \lambda_j(t+h)-h\tilde\lambda_j'(t+h)+\frac{h^2}{2}\tilde\lambda_j''(t+h)+o(h^2),  
\end{eqnarray*}  
as $h\rightarrow 0$, which implies   
\[  
\lim_{h\rightarrow 0}\frac{1}{h}(\tilde \lambda'_j(t+h)-\tilde \lambda'_j(t)) = \lim_{h\rightarrow 0}\frac{1}{2}(\tilde\lambda''_j(t+h)+\tilde\lambda''_j(t)) 
=\tilde \lambda''_j(t) 
\] 
for any $t\in I$. Hence these $\tilde C_j$'s are $C^2$.
\end{proof} 

\vspace{5mm} 

{\bf Acknowledgements.}  The authors thank the referees for valuable
comments on the original manuscprit.

%\section*{References}

%%% Please do NOT use ``\bysame'' command in your bibliography list.

\end{document}